\documentclass[a4paper]{ifacconf}
\usepackage{graphicx,amsmath,url}      
\usepackage[round]{natbib}
\graphicspath{{photo/}}
\newtheorem{assumption}{Assumption} 
\newtheorem{definition}{Definition} 
\newtheorem{remark}{Remark} 
\newtheorem{corollary}{Corollary} 
\newtheorem{Lemma}{Lemma} 
\newtheorem{theorem}{Theorem} 

\usepackage{amsmath,amssymb,amsfonts}
\usepackage[ruled,Algorithm]{algorithm}
\usepackage{natbib}
\usepackage{xcolor}

\usepackage{algorithmic}[]   

\def\portugues{1} 

\ifx\portugues\undefined
\def\portugues{0}
\fi

\if\portugues0
   \usepackage[english]{babel}
  \else
   \usepackage[spanish,brazil,english]{babel}
\fi

\usepackage[T1]{fontenc}

\usepackage[utf8]{inputenc}

\usepackage{ae}

\if\portugues1

 { }
 { }
 { }
 
\fi

\begin{document}

\if\portugues1

%
\selectlanguage{english}
	
\begin{frontmatter}

\title{Compressed Momentum-based Single-Point Zeroth-Order Algorithm for Stochastic Distributed Nonconvex Optimization} 


\author[First]{Linjing Chen} 
\author[First]{Antai Xie} 
\author[Third]{Xinlei Yi}
\author[First]{Xiaoqiang Ren}
\author[First]{Xiaofan Wang}

\address[First]{School of Mechatronic Engineering and Automation, Shanghai
University, Shanghai, China. Emails:\\\{1525417062, xatai, xqren, xfwang,\}@shu.edu.cn}
\address[Third]{Department of Control Science and Engineering, College of
Electronics and Information Engineering, Tongji University, Shang-
hai, China. Email: xinleiyi@tongji.edu.cn}

\selectlanguage{english}
\renewcommand{\abstractname}{{\bf Abstract:~}}
\begin{abstract}          
This paper proposes a compressed momentum-based single-point zeroth-order algorithm for stochastic distributed nonconvex optimization, aiming to alleviate communication overhead and address the unavailability of explicit gradient information. In the proposed algorithm, each agent has access only to stochastic zeroth-order information of its local objective function, performs local stochastic updates with momentum, and exchanges compressed updates with its neighbors. We theoretically prove that, with fixed step sizes and diminishing smoothing radius, the proposed algorithm achieves the convergence rate $\mathcal{O} (\frac{1}{\sqrt[4]{T}})$ to the stationary point. With fixed step sizes and smoothing radius, it attains a faster convergence rate ${\cal O}(\frac{1}{{\sqrt {T} }})$ towards a neighborhood of the stationary point. Numerical experiments validate the effectiveness and communication efficiency of the proposed algorithm.

\vskip 1mm
\end{abstract}

\selectlanguage{english}


\selectlanguage{brazil}

\end{frontmatter}
\else
%

\begin{frontmatter}

\title{Style for SBA Conferences \& Symposia: Use Title Case for
  Paper Title\thanksref{footnoteinfo}} 

\thanks[footnoteinfo]{Sponsor and financial support acknowledgment
goes here. Paper titles should be written in uppercase and lowercase
letters, not all uppercase.}

\author[First]{First A. Author} 
\author[Second]{Second B. Author, Jr.} 
\author[Third]{Third C. Author}

\address[First]{Faculdade de Engenharia Elétrica, Universidade do Triângulo, MG, (e-mail: autor1@faceg@univt.br).}
\address[Second]{Faculdade de Engenharia de Controle \& Automação, Universidade do Futuro, RJ (e-mail: autor2@feca.unifutu.rj)}
\address[Third]{Electrical Engineering Department, 
   Seoul National University, Seoul, Korea, (e-mail: author3@snu.ac.kr)}
   
\renewcommand{\abstractname}{{\bf Abstract:~}}   
   
\begin{abstract}                
These instructions give you guidelines for preparing papers for IFAC
technical meetings. Please use this document as a template to prepare
your manuscript. For submission guidelines, follow instructions on
paper submission system as well as the event website.
\end{abstract}

\end{frontmatter}
\fi
\section{Introduction}
With the rapid growth of datasets and the expansion of network scales, distributed optimization has attracted significant attention due to its enhanced parallel computational capabilities and inherent robustness, where agents try to cooperatively solve a problem with information exchange only limited to immediate neighbors in the network. In many networked systems, such as power systems~\citep{mao2019finite} and wireless sensor networks~{\citep{kuthadi2022optimized},~\citep{lei2022distributed}}, each computational node is typically limited to accessing information from its immediate neighbors, driven by requirements for security and scalability. 

To date, a wide variety of distributed algorithms have been proposed, including distributed gradient tracking algorithms~\citep{lei2022distributed},~\citep{chen2024compressed}, distributed ADMM~\citep{carnevale2025admm},
~\citep{gholami2023admm}, 
distributed pull-push gradient-based algorithm~\citep{nedic2025ab}, \citep{dimlioglu2025communication} and so on. Additionally, a popular SGD heuristic that has proven to be crucial in many applications is the use of momentum, i.e., the use of a weighted average of past gradients instead of the current gradient~\citep{sutskever2013importance}. Subsequently,~\cite{xia2024momentum} investigated a paradigm of momentum-based systems for nonconvex constrained optimization and distributed nonconvex optimization.~\cite{li2025stochastic} proposed a gradient tracking stochastic distributed optimization algorithm with adaptive momentum. 

However, in fields such as biochemistry, commercial applications and complex engineering systems, accessing gradients may pose significant challenges due to data privacy concerns, black-box constraints, or the high computational costs involved. To handle this problem, the gradient is approximated through random sampling and finite differences in zeroth-order (ZO) optimization. Existing ZO optimization algorithms can be divided into the following categories, namely, ZO with one-point estimator and ZO with two- and multiple-point estimator.~\cite{flaxman2004online} proposed the ZO algorithm with one-point estimator that queries one function value at each iteration to estimate 
the gradient. Nowadays, some research findings are available to integrate ZO optimization with momentum-based algorithms. For instance,~\cite{huang2022accelerated} studied a class of accelerated zeroth- and first-order momentum methods for both nonconvex mini-optimization and minimax-optimization.~\cite{qian2023zeroth} explored a ZO proximal stochastic recursive
momentum algorithm. 

Note that the distributed optimization algorithms require each agent to exchange information with its neighboring agents to obtain the global information. However, the communication channel often  has limited bandwidth. Accordingly, the communication efficient strategies have been proposed, including the event-triggered schemes and compression techniques. The former mainly decrease the number of transmission rounds~\citep{yang2022event}. The latter reduce the number of transmitted bits, such as via standard uniform quantizers, including standard uniform quantizers~\citep{xu2024quantized}, unbiased compressors~\citep{condat2022ef}, and sparsification methods~\citep{wangni2018gradient}.

In this paper, motivated by the well-known momentum-based optimization algorithm and by the practical limitation that querying multiple function values is often infeasible, especially in stochastic environments, we develop a Compressed Momentum-based Single-Point Zeroth-Order (CMSPZO) algorithm for stochastic distributed nonconvex optimization. The contributions are summarized
as follows.
\begin{itemize}
\item  We propose a compressed momentum-based algorithm for distributed nonconvex ZO optimization, where gradients are estimated from noisy single-point function evaluations. While \cite{singh2021squarm} studies a first-order distributed momentum-based optimization algorithm, and \cite{mhanna2023single} considers distributed ZO optimization without communication compression, our method combines momentum, ZO gradient estimation, and compressed communication in a unified framework.
\item The algorithm achieves an $\mathcal{O} (\frac{1}{\sqrt[4]{T}})$ convergence rate toward a stationary point with fixed step sizes and diminishing smoothing radius. Moreover, it further achieves a faster convergence rate ${\cal O}(\frac{1}{{\sqrt {T} }})$ towards a neighborhood of a stationary point with fixed step sizes and smoothing radius, where $T$ is the number of iterations. These results show that provable convergence remains achievable in the more challenging compressed communication distributed ZO setting.
    \item We develop a refined analysis for the norm squared of the gradient estimator by separately controlling the consensus error, the function value difference term, and the noise term, thereby removing the bounded iterate assumption used in prior work such as~\cite{mhanna2023single}.
\end{itemize}
$Notations$: Let $\mathbb{R}$ denote the set of real numbers, and let the vectors with all entries equal to 1 or 0 be denoted by $\mathbf{1}$ and $\mathbf{0}$, respectively. The set of $n$-dimensional real vectors is denoted by $\mathbb{R}^n$, and $\mathbf{I}$ means the identity matrix. For a differentiable function $f$, $\nabla f$ represents its gradient. For two vectors $a, b \in \mathbb{R}^n$, $\langle a, b \rangle$ is the standard inner product. For a positive semi-definite matrix $\mathbf{W}$, $\lambda(\mathbf{W})$ denotes its spectral radius. 
Additionally, $\operatorname{sgn}(\cdot)$, $\circ$, $\lfloor \cdot \rfloor$, $|\cdot|$, $\mathbb{E}_X[\cdot]$ and $\mathbb{E}[\cdot]$
denote the element-wise sign, Hadamard product, floor function, absolute value, the expectation with respect to random variable $X$ and all possible random variables, respectively. $\|\cdot\|$ and $\|\cdot\|_F$ symbolize the $\ell_2$-norm and the Frobenius norm, respectively.
\section{Preliminaries and problem formulation}
\subsection{Problem Formulation}
Consider a set of agents ${\cal V} = \{ 1,2,...,n\}$ connected by a communication network. Each agent $i$ is associated with a local objective function $F_i:{\mathbb{R}^d} \mapsto \mathbb{R}$. All agents aim to collaboratively minimize the global objective function:
\begin{equation}
\mathop {\min }\limits_{x \in {\mathbb R^d}} {f}(x): = \frac{1}{n}\sum\limits_{i = 1}^n {{F_i}(x)},
\label{Section2.1-1}
\end{equation}
where 
\begin{equation}
{F_i}(x) \buildrel \Delta \over = {\mathbb{E}_{\xi_i  \sim {\cal{D}}_i}}[{f_i}(x,{\xi_i})],
\label{Section2.1-2}
\end{equation}
with $x \in {\mathbb R^d}$ denotes the model parameter and ${\xi_i}$ is the local data that follows the local distribution ${\cal D}_i$, and ${f_i}(x,{\xi_i})$ is the stochastic local cost function of each agent $i$. $\mathbb{E}_{\xi_i \sim \mathcal{D}_i}[\cdot]$ denotes the expectation with respect to the random sample $\xi_i$ drawn from $\mathcal{D}_i$ of client $i$.
\begin{assumption}\label{lplplp}
Assume that $\mathbb{E}_{\xi_i  \sim {\cal{D}}_i}{{{[{{f_i}({x_{}},{\xi _{i}})}]}^2}}\le \gamma _1$, where $\gamma _1>0$.   
\end{assumption}

\begin{assumption}\label{section 2.1-ass:the local function}
We assume the Lipschitz continuity of the local objective functions 
$f_i(x,\xi_i)$ such that
\begin{equation}
    | f_i(x,\xi_i) - f_i(x',\xi_i) | \le  L_{f_1} \| x - x' \|,
    \quad \forall x, x' \in \mathbb{R}^d,
\end{equation}
where $L_{f_1}$ denotes the Lipschitz constant.
\end{assumption}
\begin{assumption}\label{section 2.1-ass:the objective function}
Both $\nabla F_i(x)$ and $\nabla^2 F_i(x)$ exist and are continuous, and there exists a constant $L_{{f_2}} > 0$ such that
$\|\nabla^2 F_i(x)\| \leq L_{{f_2}}$, $\forall x \in \mathbb{R}^d.$
\end{assumption}
From \textit{Assumption~\ref{section 2.1-ass:the objective function}}, we can say that the objective function 
$f(x)$ is $L_{{f_2}}$-smooth:
\begin{equation}
f(x) \le f(x') + \langle \nabla f(x'), x - x' \rangle + \frac{L_{{f_2}}}{2} \| x - x' \|^2,
\end{equation}
or alternatively,
\begin{equation}
\|\nabla f(x)-\nabla f(x')\| \le {L_{{f_2}}}\|x-x'\|,\quad \forall x,x' \in \mathbb{R}^d.
\end{equation}
To address the problem \eqref{Section2.1-1}, we focus on the following aspects: exchanging information with other agents and collecting information about the local cost function. We introduce the following mechanisms to address these two aspects, together with several additional assumptions.
\subsection{Graph Theory}
The network is described by an
undirected and connected graph ${\cal G} = ({\cal V},{\cal E})$, meaning communication
links work in both directions, and between any two agents,
we can find a path of links. 
Define ${\cal E} \subset {\cal V} \times {\cal V}$ as the set of communication channels. An channel $(i,j) \in {\cal E}$ indicates that the agent $i$ can communicate information with the agent $j$. The non-negative weighted matrix $\mathbf{W}$ is represented as $\mathbf{W} = [{w_{ij}}] \in {{\mathbb{R}}^{n \times n}}$. ${w_{ij}} > 0$ if $(i,j) \in {\cal E}$, and otherwise ${w_{ij}} = 0$. And all diagonal elements $w_{ii}$ are strictly positive. The set of neighbors of the agent is denoted by ${{\cal N}_i} = \{ j|(i,j) \in {\cal E}\}$. 
\begin{assumption}\label{ass:graph}
 The weighted matrix $\mathbf{W}$ is symmetric and doubly stochastic, i.e., ${\mathbf{1}^T}\mathbf{W} = {\mathbf{1}^T},\mathbf{W}\mathbf{1} = \mathbf{1}$. 
\end{assumption}
It is well known (see \citep{koloskova2019decentralized}) that for the matrix $\mathbf{W}$ associated with ${\cal G}$, its eigenvalues satisfy 
\begin{equation*}
1=|\lambda_1(\mathbf{W})|>|\lambda_2(\mathbf{W})|\ge\cdots\ge|\lambda_n(\mathbf{W})|,
\end{equation*}
and the spectral gap is defined as
$\delta= 1 -|\lambda_2(\mathbf{W})| \in (0,1]$, where $\lambda_1(\mathbf{W})$ is the eigenvalue of $\mathbf{W}$ with the largest magnitude, and $\{\lambda_i(\mathbf{W})\}_{i=2}^n$ are the remaining eigenvalues.
\subsection{Zeroth-Order Optimization}
This section concentrates on utilizing the zeroth-order
(ZO) information to estimate the gradients of local objective
functions. In the absence of explicit gradient expressions, we assume each agent could query the local function value at one point.

Let $F(x) = \mathbb{E}_{\xi}[f(x,\xi)]$, where $f(x,\xi)$ denotes a stochastic differentiable function.  As shown in \cite{flaxman2004online}, the one-point of the gradient can be constructed as
\begin{equation}
g(x,\xi): = \frac{d}{\gamma }f(x + \gamma u,\xi)u,\label{section2.3-one point in Flaxman}
\end{equation}
where $\gamma  > 0$ denotes the smoothing radius and $u$ is a random vector that follows a symmetric distribution. 

To solve the problem (\ref{Section2.1-1}), agent $i$ can only query the function value of $f_i$ at one point. The function query is assumed to be noisy $\tilde{f}_i({x_{i}},{\xi _{i}})= f_i({x_{i}},{\xi _{i}})+ \varphi_{i}$, where $\varphi_i$ denotes additive observation noise. 
\begin{remark}\label{remark 1}
Since ZO optimization methods require a function evaluation at each iteration to estimate
the gradient, the queried function value may be affected not only by the stochastic nature of $f_i$, but also an additional noise introduced in the function query process. Thus, $\varphi_i$ is introduced to model this additive observation noise. This separation clarifies how different noise sources affect the gradient estimator and the convergence analysis.
\end{remark}
Then, we assume that the function queries are 
\begin{equation*}
    \tilde{f}_i ({x_{i,t}} + {\gamma _g}{u_{i,t}},{\xi _{i,t}})= f_i ({x_{i,t}} + {\gamma _g}{u_{i,t}},{\xi _{i,t}})+ \varphi_{i,t},
\end{equation*}
where $x_{i,t}\in \mathbb{R}^d$ is the estimation of $x$ maintained by agent $i$ at time step $t$, $\varphi_{i,t}$ denotes the additive noise at time step $t$, $\gamma _g>0$ is the smoothing radius, $u_{i,t}$ is a random perturbation vector at time step $t$. 

Inspired by (\ref{section2.3-one point in Flaxman}), we have
\begin{align} 
g_{i,t} &=\frac{d}{\gamma_g}  \tilde{f}_i \big( x_{i,t} + \gamma_g u_{i,t}, \xi _{i,t} \big)u_{i,t}\nonumber\\
&=\frac{d}{\gamma_g} (f_i ( x_{i,t} + \gamma_g u_{i,t}, \xi _{i,t} ) + \varphi _{i,t} )u_{i,t}.\label{section2.3-one point in this paper}
\end{align}

For the noise term $\varphi_{i,t}$ and the perturbation vector $u_{i,t}$, we adopt the following assumptions. 
\begin{assumption}\label{section2.3-ass:noise}
The noise term $\varphi_{i,t}$ is a zero-mean un-
correlated noise with bounded variance, meaning $\mathbb{E}[\varphi_{i,t}]=0$ and $\mathbb{E}[\varphi_{i,t}^2] = \vartheta_1 < \infty, i \in \mathcal{V}$, and $\mathbb{E}[\varphi_{i,t}\varphi_{j,t}] = 0$ for $i\ne j$.
\end{assumption}
\begin{remark}
\textit{Assumption \ref{section2.3-ass:noise}} is common in the literature (see \citep{bychkov2024accelerated}). $\varphi_i$ is modeled as additive, zero-mean noise with bounded second moment and is assumed to be uncorrelated across agents. While the expectation of the local function might not change with this addition, its variance increases, presenting an added difficulty to the analysis.
\end{remark}
\begin{assumption}\label{section2.3-ass:perturbation vector}
The perturbation vector $u_{i,t} = (u_{i,t,1}, \ldots,\\
u_{i,t,d})^\top \in {\mathbb R^d}$ is chosen independently by each agent $ i \in \mathcal{N} $ from others and previous samples. In addition, the elements of $u_{i,t}$ are assumed i.i.d. with
$\mathbb{E}[u_{i,t,d}^2] = \sigma_1 > 0$, 
which implies that $\mathbb{E}[u_{i,t}u_{i,t}^T] = \sigma_1 I_d$. We assume that there exists a constant $\sigma_2 > 0$ such that $\|u_{i,t}\| \le \sigma_2.$
\end{assumption}
\begin{remark}
An example of a perturbation vector satisfying \textit{Assumption \ref{section2.3-ass:perturbation vector}} is to choose every dimension of $u_{i,t}$ from the symmetric Bernoulli distribution on $\{-\tfrac{1}{\sqrt{d}}, \tfrac{1}{\sqrt{d}}\}$. Then, $\sigma_1=\frac{1}{d}, \sigma_2=1$.
\end{remark}
Before proceeding with the analysis, we define ${{\cal L}_t}$ as the $\sigma$-algebra generated by 
\begin{equation*}
{\{\{u_{i,0},\xi_{i,0},\varphi_{i,0}\}_{i=1}^n,\,
\dots,\,
\{u_{i,t-1},\xi_{i,t-1},\varphi_{i,t-1}\}_{i=1}^n\}}.
\end{equation*}
By construction, $\{\mathcal{L}_t\}_{t\ge 0}$ forms an increasing filtration, i.e., $\mathcal{L}_{t-1} \subseteq \mathcal{L}_t$ for all $t \ge 1$.
\begin{Lemma} (Proposition 3.3. of \citep{mhanna2023single})\label{section2.3-lemma1}
Under \textit{Assumptions \ref{section2.3-ass:noise}--\ref{section2.3-ass:perturbation vector}}, then $g_{i,t}$ is a biased estimator of the agent's gradient $\nabla {F_i}({x_{i,t}}),\forall i \in {\cal V}$ for every $t \ge  0$, i.e., 
\begin{equation}\label{section2.3-lemma1-1}
\mathbb{E}_{{{\cal L}_t}}[{g_{i,t}}] = d{\sigma _1}(\nabla {F_i}({x_{i,t}}) + {b_{i,t}})
\end{equation}
with
\begin{equation}\label{section2.3-lemma1-2}
\begin{aligned}
{{b_{i,t}}}={\frac{{{\gamma _g}}}{{2{\sigma _1}}} \mathbb{E}{_u}[{u_{i,t}}u_{i,t}^T{\nabla ^2}{F_i}(\upsilon){u_{i,t}}}],
\end{aligned}
\end{equation}
where ${{b_{i,t}}} \in \mathbb{R}^d$ is the bias with respect to the true gradient, and 
$\upsilon  \in [x_{i,t},\, x_{i,t} + \gamma_g u_{i,t}]$. 
\end{Lemma}
\subsection{Compression Operators}
To reduce communication overhead, we consider a scenario where information exchange between agents is compressed. Specifically, we consider the  compressors $\mathcal{C}(\cdot)$ with bounded relative compression error satisfying the following assumption.  
\begin{definition}\label{section2.4-definition-1}
A function $\mathcal{C}: \mathbb{R}^d \to \mathbb{R}^d$ is called a compression operator, if there exists a positive constant $\omega \in (0,1]$, such that for every $x \in \mathbb{R}^d$:
\begin{equation}
\mathbb{E}_{\mathcal{C}}[\|x - \mathcal{C}(x)\|_2^2] 
\le (1-\omega)\|x\|_2^2,
\end{equation}
where expectation is taken over the randomness of $\mathcal{C}$. 
Without loss of generality, we assume that $\mathcal{C}(0)=0$.
\begin{remark}
As stated in~\cite{liao2022compressed} is naturally satisfied by various compressors in distributed algorithms. For example, the Rand$_k$ sparsifier~\citep{beznosikov2023biased}, the
Top$_k$ sparsifier~\citep{zou2022downlink} and norm-sign compressor~\citep{yi2022communication}.
\end{remark}
\end{definition}
\section{Main results}
In this section, we propose a Compressed Momentum-based Single-Point Zeroth-Order (CMSPZO) algorithm. First, we introduce the CMSPZO algorithm and its motivation in Section~\ref{section3.1}. Then, we provide the supporting lemmas and establish the convergence in Section~\ref{section3.2}. 
\subsection{Algorithm Description}\label{section3.1}
When the gradient information is available, \cite{yu2019linear} proposed the following first-order momentum-based algorithm 
\begin{subequations}
    \begin{align}
m_{i,t} &= \beta m_{i,{t-1}} + 
\nabla F_i(x_{i,t-1}, \xi_{i,t-1}),\label{section3.1-momentum-based algorithm in other paper-2}\\
x_{i,{t}}&=\sum_{j=1}^{n}\omega_{ij}(x_{j,{t-1}}-\eta m_{j,t}),
\end{align}
\end{subequations}
where $m_{i,t}$ is the local momentum of agent $i$ at time step $t$, $\eta$ is the step size, $\beta \in [0,1)$ is momentum coefficient. 

However, first-order information may not always be available in practice, as gradient computation can be impractical or too costly for large-scale or highly complex systems. Then, we
replace the gradient term in (\ref{section3.1-momentum-based algorithm in other paper-2}) with the single-point ZO gradient estimator in~(\ref{section2.3-one point in this paper}). 

To realize the exchange of compressed parameters between agents, each agent $i$ maintains an estimate $\hat{x}_{i,t}$ of $x_{i,t}$. 
We then propose the CMSPZO algorithm (Algorithm 1), which integrates the single-point ZO estimator and momentum-based method. 

At iteration $t$, each agent $i$ computes a single-point ZO gradient estimator $g_{i,t}$ and performs a local momentum-based update to form ${x_{i,t+\frac{1}{2}}}$. It then transmits compressed information ${\cal C}(x_{i,t+\frac{1}{2}} - \hat{x}_{i,t})$ to its neighbors. Next, each agent updates its local variables according to~(\ref{section3.1-algorithm1-neighbors}). Consequently, the agents finally take a consensus step in~(\ref{section3.1-algorithm1-consensus}) with appropriate weighting decided by entries of \textbf{W}.

\begin{algorithm}[htbp]
\caption{\underline{C}ompressed \underline{M}omentum-based \underline{S}ingle-\underline{P}oint \underline{Z}eroth-\underline{O}rder (CMSPZO) Algorithm}
\label{alg:cmsgd1p}
\begin{algorithmic}[1]   
  \STATE \textbf{Input:} Stopping time $T$, adjacency matrix $\mathbf W$, and positive parameters $\gamma_x$, $\eta$, $\beta\in [0,1)$, $\gamma_g$.
  \STATE \textbf{Initialization:} For all $i \in \mathcal V$, initial variable
         $x_{i,0} \in \mathbb R^d$, ${x}_{i,-1} =\textbf{0}$, $\hat{x}_{i,-1}=\textbf{0}$, $m_{i,-1} = \textbf{0} \in \mathbb R^d$.
  \FOR{$t = 0$ to $T-1$ \textbf{ in parallel on all agents} $i$}
    \STATE One-point ZO gradient:
      \begin{itemize}
        \item Sample $\xi_{i,t}$ from the distribution $\mathcal D_i$;
        \item Sample $u_{i,t}$ from the symmetrical distribution;
        \item Query $f_i(x_{i,t} + \gamma_g u_{i,t}, \xi_{i,t})$;
        \item Compute stochastic gradient $g_{i,t}$ using (\ref{section2.3-one point in this paper}).
      \end{itemize}
    \STATE Update $m_{i,t}$ by
    \begin{equation}
   {m_{i,t}} =  \beta {m_{i,t - 1}} +{g_{i,t}}.\label{section3.1-algorithm1-estimete lacal momentum}
\end{equation}
   \STATE Update $x_{i,t+\frac{1}{2}}$ according to
   \begin{equation}
{x_{i,t+\frac{1}{2}}} = {x_{i,t}} - \eta (\beta{m_{i,t}}+{g_{i,t}}), \label{section3.1-algorithm1-local SGD Update}
\end{equation}
    \STATE Compute ${\cal C}(x_{i,t+\frac{1}{2}} - \hat{x}_{i,t})$ and broadcast it to its neighbors $\mathcal N_i$.
\STATE Receive ${\cal C}(x_{j,t+\frac{1}{2}} - \hat{x}_{j,t})$ from $j \in \mathcal N_i$.
\STATE Update ${\hat x_{j,t+ 1}}$ by
\begin{equation}
{\hat x_{j,t + 1}} = {\hat x_{j,t}} + {\cal C}({x_{j,t+\frac{1}{2}}} - {\hat x_{j,t}}).\label{section3.1-algorithm1-neighbors}
\end{equation}
    \STATE \textbf{Consensus:}
\begin{equation} 
{x_{i,t + 1}} = {x_{i,t +\frac{1}{2}}}+ {\gamma _x}\sum\limits_{j \in {{\cal N}_i}} {{w_{ij}}} ({\hat x_{j,t + 1}} - {\hat x_{i,t + 1}}). \label{section3.1-algorithm1-consensus}
\end{equation}
  \ENDFOR
  \STATE \textbf{Output:} $\{x_i(T)\}$.
\end{algorithmic}
\end{algorithm}
\subsection{Supporting Lemmas}\label{section3.2}
The following lemmas are used in the proofs.

For the sake of analysis, we introduce the following notations: 
\begin{align}\label{whh}
{{\mathbf{x}_t}} &= [{x_{1,t}^T},...,{x_{n,t}^T}]^T,~~{{\bar {\mathbf{x}}}_t} = [{{\bar x}_t^T},...,{{\bar x}_t^T}]^T, \nonumber \\
{{\hat{\mathbf{x}}_t}}&=[{{\hat x}_{1,t}^T},...,{{\hat x}_{n,t}^T}]^T,~~{\mathbf{m}_t} =[{m_{1,t}^T},...,{m_{n,t}^T}] ^T, \nonumber\\
{\mathbf{g}_t} &= [{g_{1,t}^T},...,{g_{n,t}^T}]^T,~~ 
\mathcal{C}(\mathbf{x}) = [\mathcal{C}(x_{1})^T,...,\mathcal{C}(x_{n})^T]^T,
\end{align}   
where ${\bar x}_t = \frac{1}{n}\sum\limits_{i = 1}^n{{x_{i,t}}}$.
To begin with, we define
\begin{equation}\label{section3.2-Lyapunov candidate function}
{\chi _t} = \mathbb{E}[\| {{\mathbf{x}_t} - {{\bar {\mathbf{x}}}_t}} \|_F^2] + \mathbb{E}[\| {{\mathbf{x}_t} - {{\hat {\mathbf{x}}}_t}} \|_F^2],    
\end{equation}
where incorporates the consensus error term $\mathbb{E}[\| {{\mathbf{x}_t} - {{\bar {\mathbf{x}}}_t}} \|_F^2]$ and the compressed error term $\mathbb{E}[\| {{\mathbf{x}_t} - {{\hat {\mathbf{x}}}_t}} \|_F^2]$.
They play the key role for obtaining the convergence results of Algorithm 1.   
\begin{Lemma}(Consensus error)\label{section3.2-lemma2-consensus error}
Let $\{\mathbf{x}_t\}$ be the sequence generated by Algorithm 1, we have 
\begin{align*}
 \mathbb{E}&[\| \mathbf{x}_t - \bar {\mathbf{x}}_t \|_F^2]\nonumber\\
 \le& {\varepsilon _1}(1 + \alpha _1^{ - 1})\mathbb{E}[\| {{{{\mathbf{x}}}_{t - 1}} - {\bar {\mathbf{x}}_{t - 1}}} \|_F^2]+{{{{\varepsilon _4}\mathbb{E}}}{[\| {{{\mathbf{x}}_{t - 1}} - {{\hat {\mathbf{x}}}_{t - 1}}} \|_F^2]}}\nonumber
  \end{align*}
 \begin{align}
\label{section3.2-lemma2-consensus error-1}
 &+((1 + {\alpha _1})({\varepsilon _1} + {\varepsilon _2})+{\varepsilon _3})\mathbb{E}[\| {\eta (\beta {{\mathbf{m}}_{t - 1}} + {{\mathbf{g}}_{t - 1}})} \|_F^2].
   \end{align}
where $\alpha_1>0$ and ${\varepsilon_1}$--${\varepsilon_4}$ are given in Appendix~\ref{appendix6.2-Proof of Lemmas}. 

\textbf{Proof.} See Appendix~\ref{appendix6.2-Proof of Lemmas} for the proof. 
\end{Lemma}
\begin{Lemma}(Compressed error)\label{section3.2-lemma2-compressed error}
Let $\{\mathbf{x}_t\}$ be the sequence generated by Algorithm 1, we have
\begin{equation}\label{section3.2-lemma2-compressed error-1}
\begin{aligned}
 &\mathbb{E}[\| {{{\mathbf{x}}_t} - {{\hat {\mathbf{x}}}_t}} \|_F^2]
\le {\kappa _3}\mathbb{E}[\| {{{\mathbf{x}}_{t - 1}} - {{\bar {\mathbf{x}}}_{t - 1}}} \|_F^2]\\
  &+( (1 + {\alpha _2})(1 - \omega )(1 + {\alpha _3})+{\kappa _2})\mathbb{E}[\| {{{\mathbf{x}}_{t - 1}} - {{\hat {\mathbf{x}}}_{t - 1}}} \|_F^2]\\
&+((1 + {\alpha _2})(1 - \omega )(1 + \alpha _3^{ - 1})+ {\kappa _1})\mathbb{E}[\| { \eta (\beta {\mathbf{m}_{t-1}} + {\mathbf{g}_{t-1}})} \|_F^2].
\end{aligned}
\end{equation}
where $\alpha_2, \alpha_3>0$ and ${\kappa_1}$--${\kappa_3}$ are provided in Appendix~\ref{appendix6.2-Proof of Lemmas}. 

\textbf{Proof.} See Appendix~\ref{appendix6.2-Proof of Lemmas} for the proof. 
\end{Lemma}
\begin{Lemma}\label{section3.2-lemma bound of gradient for lemma 2-3}
Suppose \textit{Assumptions \ref{section2.3-ass:noise}--\ref{section2.3-ass:perturbation vector}} hold and $\textbf{m}_{-1}=\textbf{0}$, we have
\begin{align}\label{section3.2-lemma4-whole}
\mathbb{E}&[\| {\eta (\beta {\mathbf{m}_{t - 1}} + {\mathbf{g}_{t - 1}})} \|_F^2] \le  \frac{{16{\eta ^2}{\beta ^2}{d^2}\sigma _2^2L_{{f_1}}^2}}{{(1 - \beta )\gamma _g^2}}\sum\limits_{k = 0}^{t - 1} {{\beta ^{t - 1 - k}}} {\chi _k}\nonumber\\
&+ \frac{{16{\eta ^2}{d^2}\sigma _2^2L_{{f_1}}^2}}{{\gamma _g^2}}{\chi _{t-1}} + 2{\eta ^2}{\kappa _4}\left(1 + \frac{{{\beta ^2}}}{{{{(1 - \beta )}^2}}}\right),  
\end{align}
\end{Lemma}
where ${\kappa _4}=\frac{{4n{d^2}\sigma _2^2}(L_{{f_1}}^2\gamma _g^2\sigma _2^2 + 2\gamma _1 + 2{\vartheta _1})}{{\gamma _g^2}}$.

\textbf{Proof.} See Appendix~\ref{appendix6.2-Proof of Lemmas} for the proof. 

Now, we show the convergence of CMSPZO algorithm under general nonconvex settings. 
\begin{theorem}
 Suppose \textit{Assumptions \ref{lplplp}--\ref{section2.3-ass:perturbation vector}} hold and $f^*=\inf_{x \in \mathbb{R}^d} f(x)>-\infty$. For the total number of iterations $T >{\tilde \varepsilon _1}$, each agent runs Algorithm \ref{alg:cmsgd1p} 
with 
\begin{align*}
     {m_2}<{\gamma _x} < {m_1},~\eta  = (1 - \beta )\sqrt {\frac{n}{T}}, 
 \end{align*} 
where
 \begin{align}
 &{{m_1}} =~{\left(\frac{{3\phi {\omega ^3}}}{{2({\delta ^3}{\omega ^3} + 500\delta{\lambda ^2} + {{\tilde \varepsilon }_2})}}\right)^{\frac{1}{3}}},\nonumber\\
 &{{m_2}}=\frac{{{{\tilde \varepsilon }_4} + \sqrt {\tilde \varepsilon _4^2 + \frac{{27( - {{\tilde \varepsilon }_3})\phi }}{2\delta }} }}{{2( - {{\tilde \varepsilon }_3})}},\nonumber\\
 &\rho= \frac{{16{\eta ^2}{d^2}\sigma _2^2L_{{f_1}}^2}}{{\gamma _g^2}},~~~\phi = \frac{{2\rho }}{{{{(1 - \beta )}^2}}},\nonumber
\end{align}
and ${\tilde \varepsilon }_1$--${\tilde \varepsilon }_4$
are provided in Appendix~6.4. 
Then, \begin{align}
  \frac{1}{T}&\sum\limits_{t = 0}^{T - 1} {\sum\limits_{i = 1}^n {\mathbb{E}{{[\| {{x_{i,t}} - {{\bar x}_t}} \|}^2]}} }={\cal {O}}(\frac{1}{{T}})+ {\cal{O}}(\frac{1}{{T}\gamma _g^2}), \label{section3.2-theorem1-consensus error}\\
 \frac{1}{T}&\sum\limits_{t = 0}^{T - 1} \mathbb{E}{{[\| {\nabla {f}({{\bar x}_t})} \|}^2]}={\cal O}({{\gamma _g^2}})+ {\cal O}(\frac{1}{{\sqrt {nT} }}) + {\cal O}(d_1\sqrt {\frac{n}{T}})\nonumber \\
 &+ {\cal O}(\frac{d_2}{\gamma _g^2}\sqrt {\frac{n}{T}}),\label{section3.2-theorem1-gradient square norm}
 \end{align}
 \end{theorem}
 where $d_1$--$d_2$ are given in Appendix~6.4.

\textbf{Proof.} The detailed proof is given in Appendix~6.4.
\begin{remark}
It is worth noting that, in contrast to~\cite{mhanna2023single}, we analyze the norm squared of the gradient estimator without assuming that $\|\mathbf{x}_t\|<\infty$ almost surely.
\begin{remark}
Note that, $\gamma_g$ controls a bias-variance trade-off in \eqref{section3.2-theorem1-gradient square norm}. When it is fixed, the bound in \eqref{section3.2-theorem1-gradient square norm} is governed by the terms $\mathcal{O}(1/\sqrt{nT})$ and $\mathcal{O}(\sqrt{n/T})$, while the term ${\cal O}({{\gamma _g^2}})$ characterizes the smoothing bias. And a large $\gamma_g$ increases the bias term ${\cal O}({{\gamma _g^2}})$, whereas choosing $\gamma_g$ too small enlarges the last terms in (\ref{section3.2-theorem1-consensus error})--\eqref{section3.2-theorem1-gradient square norm}. 
\end{remark}
\end{remark}
From the right-hand side of (\ref{section3.2-theorem1-gradient square norm}), it follows that sublinear convergence rate can be fulfilled if ${\gamma _g}$ is chosen as a diminishing smoothing radius, that is, ${\gamma _g} ={T^{ - \frac{1}{8}}}$, which is presented in the following result.
\begin{corollary}
Under the same assumptions and parameters settings in Theorem 1, let ${\gamma _g} ={T^{ - \frac{1}{8}}}$ and $T > {\tilde \varepsilon _5}$, then, 
\begin{align}
 \frac{1}{T}&\sum\limits_{t = 0}^{T - 1} {\sum\limits_{i = 1}^n {\mathbb{E}{{[\| {{x_{i,t}} - {{\bar x}_t}} \|}^2]}} }= {\cal{O}}(\frac{1}{{T}})+{\cal O}(\frac{1}{{\sqrt[4]{{{T^3}}}}}),\label{section3.2-theorem1-consensus error corollary} \\ 
\frac{1}{T}&\sum\limits_{t = 0}^{T - 1} \mathbb{E}{{[\| {\nabla {f}({{\bar x}_t})} \|}^2]}={\cal O}(\frac{1}{{\sqrt {nT} }})  + {\cal O}(d_1\sqrt {\frac{n}{T}}) + {\cal O}(\frac{{1}}{{\sqrt[4]{T}}})\nonumber\\
&+ {\cal O}(\frac{{d_2\sqrt {{n}} }}{{\sqrt[4]{{{T}}}}}),\label{section3.2-theorem1-corollary1-gradient square norm}
\end{align}
\end{corollary}
where ${\tilde \varepsilon _5}$ can be found in Appendix~6.5.

\textbf{Proof}. The detailed proof is given in Appendix~\ref{1234567}.
\section{Simulations}
In this section, we provide the results of the simulations to prove the effectiveness of our algorithm. 

First, the system consists of $n=6$ agents and dimension $d = 30$ that exchange information through a connected undirected communication graph. And our theoretical results via numerical experiments on a nonconvex distributed binary classification task~\citep{xie2024communication}, \citep{wang2025compressed}:
\begin{equation}
f_i(x, \xi_i)
= \frac{n}{{{m_i}}}\sum_{j=1}^{m_i} 
\log(1 + e^{-p_{ij} x^\top q_{ij}})+\sum_{\iota=1}^{n} 
\frac{\varpi \kappa  [x]_\iota^2}{1 + \kappa [x]_\iota^2},
\end{equation}
where the first term
is the logistic regression function, $m_i = 200$ denotes the number of local observations per agent. 
$\xi_i = (p_{ij}, q_{ij})$, 
$p_{ij} \in \mathbb{R}^d$ is the Gaussian features regenerated per iteration, $q_{ij} \in \{-1, 1\}$ is labels. 
The second term represents a nonconvex regularizer with the regularization parameters $\varpi = 0.001$ and $ \kappa = 1$. $[x]_\iota$ is the $\iota$-th coordinate of $x \in \mathbb{R}^d$. 

Moreover, we adopt the unbiased 2-bit quantizer in \cite{wang2025compressed}, which satisfies \textit{Definition \ref{section2.4-definition-1}}. The quantization operator 
$\mathcal{C}(\cdot)$ takes the form
\begin{equation}
    \mathcal{C}(x)
    = \frac{\|x\|_{\infty}}{2^{{k_2}-1}} \mathrm{sgn}(x) \circ
      \left \lfloor 
      \frac{2^{{{k_2}}-1}|x|}{\|x\|_{\infty}} +{\varpi_1}\right \rfloor,
\end{equation}
where ${k_2}=2$, $\varpi_1$ is the random perturbation vector uniformly chosen from $[0,1]^d$.

The performance metric 
$P(T) = \min_{t\in T} \{\mathbb{E}_{\xi}[\|\nabla f(\bar{x}_t,\xi)\|^2]
+ \tfrac{1}{n}\sum_{i=1}^n \|x_{i,t} - \bar{x}_t\|^2 \}$ is evaluated in terms of both the iteration rounds and the communication bits among agents. The hyperparameters used in the experiments are summarized in Table 1. 

We compare the Algorithm 1 against distributed stochastic
gradient tracking method with
a single-point gradient estimator (DSGT-1P) in \citealp{mhanna2023single} and single-point zeroth-order stochastic compressed distributed primal-dual algorithm (ZSC-PD). 
\begin{table}
\centering
\caption{Hyperparameters for Algorithm 1}
\label{tab:hyper}
\begin{tabular}{c c c c c c}
\hline
Algorithm & $\gamma_x$ & $\eta$ & $\beta$& $\gamma_g$ \\
\hline
CMSPZO & $0.3$ & $0.02$ & $0.9$ & $0.99$ \\
ZSC-PD & $0.3$ &0.02& $-$ & $0.99$  \\
DSGT-1P & $-$ & $0.02$ & $-$ & $0.99$ \\
\hline
\end{tabular}
\end{table}
\begin{figure}
\begin{center}
\includegraphics[width=3.5in,height=2.3in]{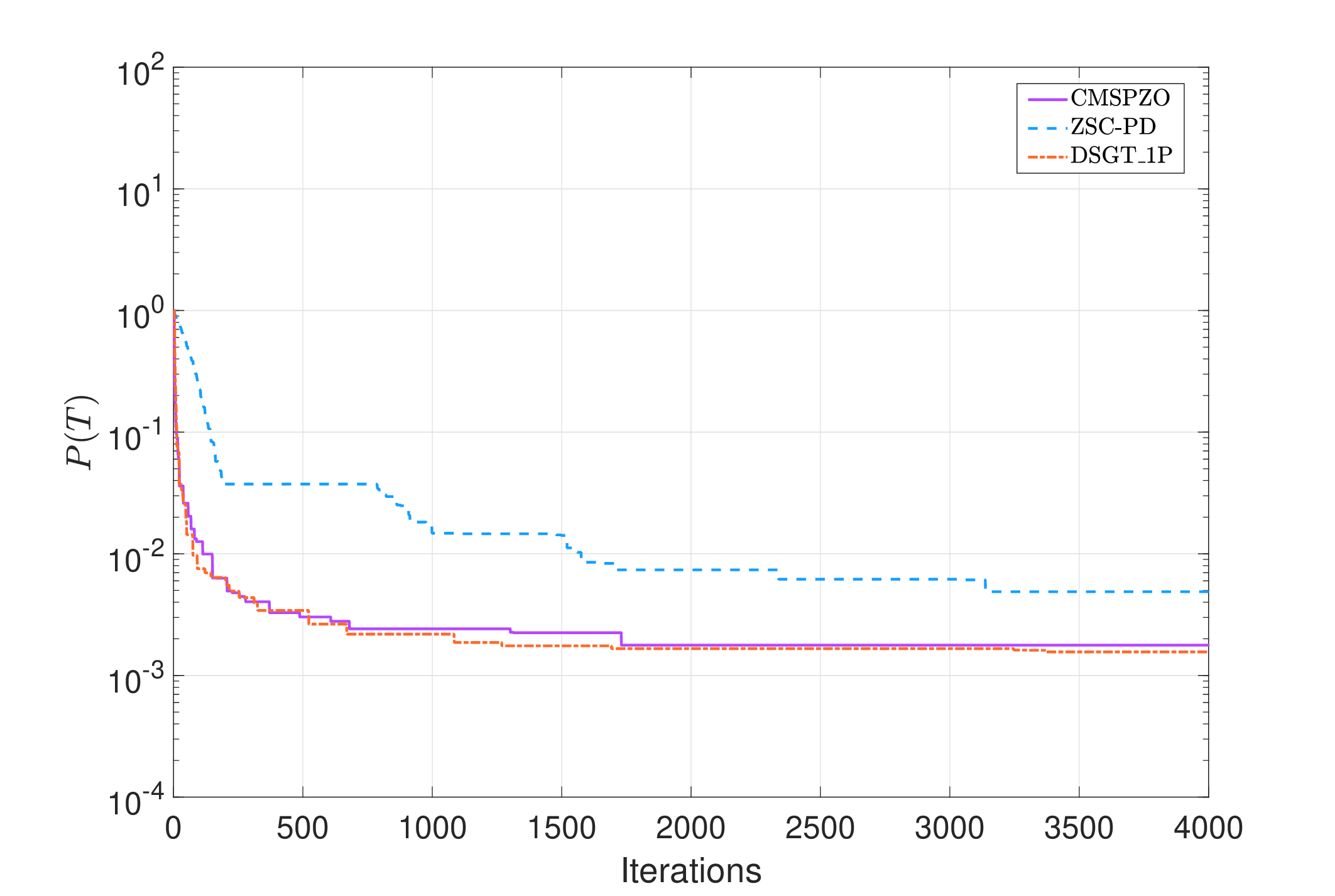}
\caption{Evolutions of $P(T)$ with respect to the number of iterations.}
\end{center}
\end{figure}
\begin{figure}
\begin{center}
\includegraphics[width=3.5in,height=2.3in]{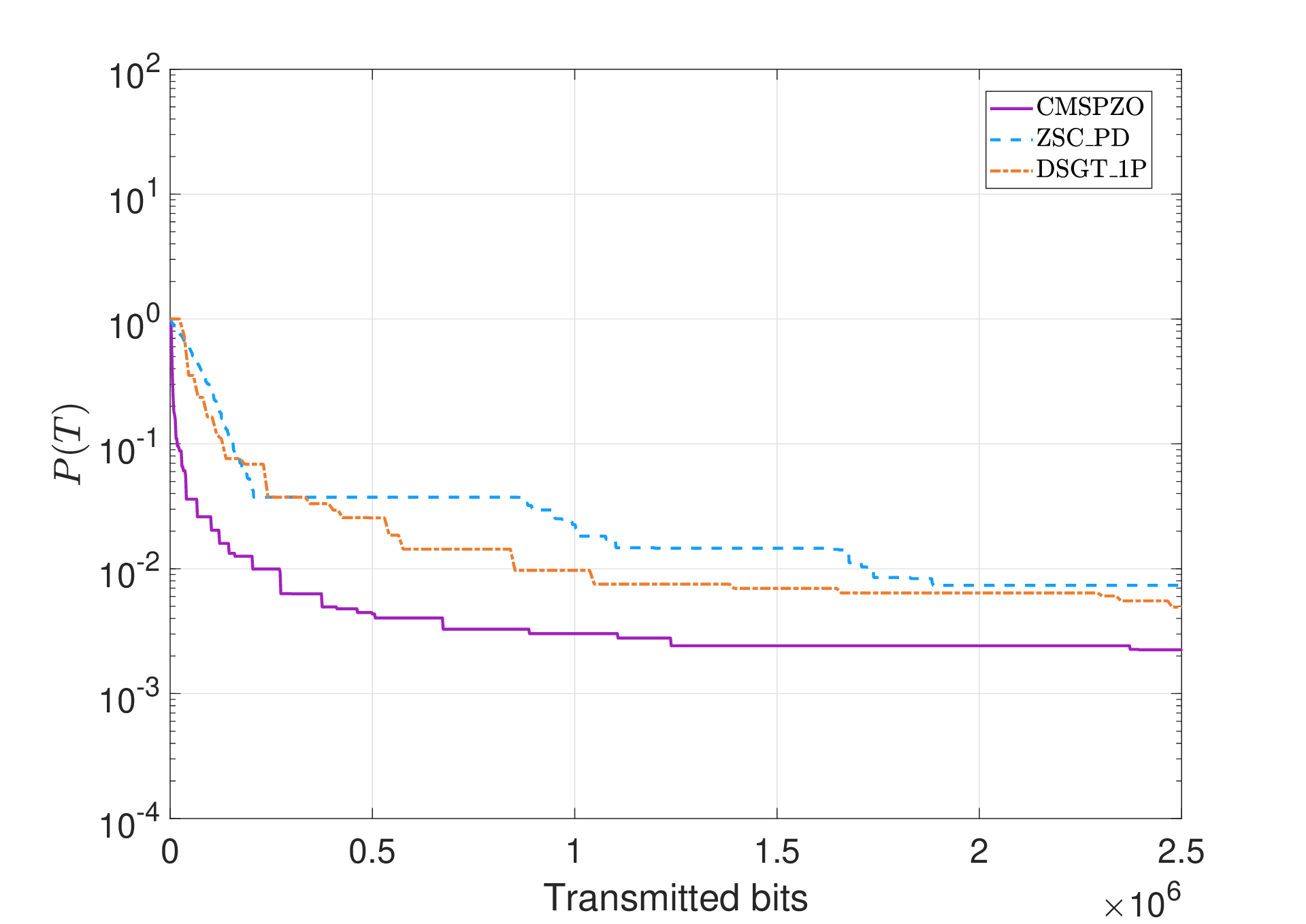}
\caption{Evolutions of $P(T)$ with respect to the number of transmitted bits.}
\end{center}
\end{figure}
 \begin{figure}
 \begin{center} \includegraphics[width=3.5in,height=2.6in]{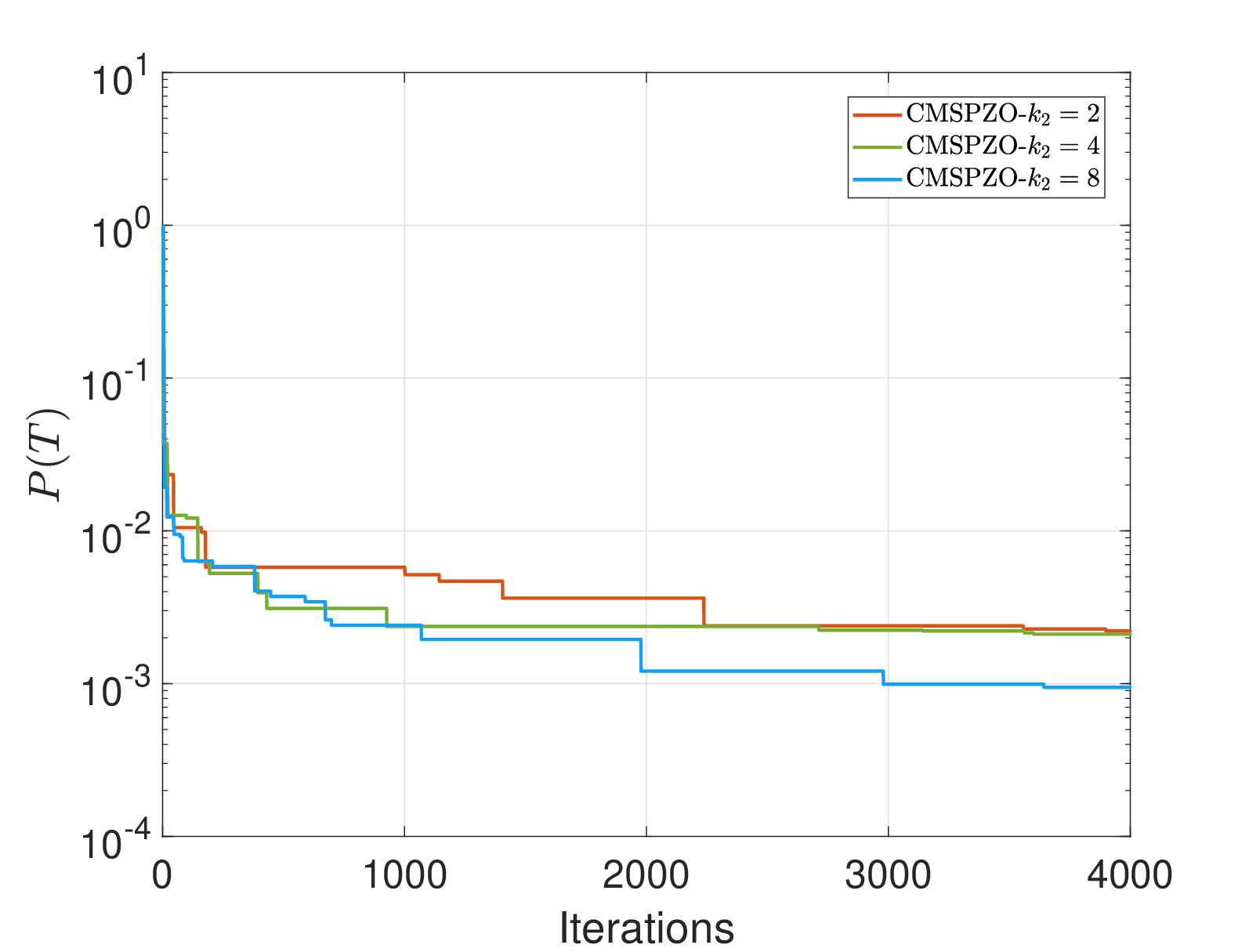}
 \caption{Evolutions of $P(T)$ with respect to $k_2$.}
 \end{center}
 \end{figure}
 
 Fig. 1 shows that CMSPZO converges slightly slower than DSGT-1P but still faster than ZSC-PD. Furthermore, Fig. 2 illustrates that compared to other algorithms, CMSPZO converges to the same accuracy with fewer bits. This also demonstrates the efficiency of CMSPZO algorithm. Fig. 3 shows that the relationship between the quantization bits $k_2$ and the convergence rates. As quantization bits
increases, i.e., more transmitted data, the convergence speed becomes faster.
\section{Conclusion}
In this paper, we proposed the CMSPZO algorithm for stochastic distributed nonconvex optimization under the assumption that the gradient is not available and that only noisy single queries of the objective function are available at a time. We confirmed that the proposed algorithm achieves the convergence rate $\mathcal{O} (\frac{1}{\sqrt[4]{T}})$ with the exact solution under fixed step sizes and diminishing smoothing radius. And it can achieve the sublinear convergence rate ${\cal O}(\frac{1}{{\sqrt {T} }})$ towards a neighborhood of the stationary point under fixed step sizes and smoothing radius. 
Future directions
include extending the algorithm to the one-point residual-feedback oracle.
\bibliography{ifacconf} 
\section{APPENDIX} 
\subsection{Useful Lemmas}\label{appendix6.1-Useful Lemmas}
The following results are used in the proofs. More notably, unless specified otherwise, for a vector $\mathbf{u}$, we write $\|\mathbf{u}\|$ to denote the $\ell_2$-norm $\|\mathbf{u}\|_2$.
\begin{Lemma} \citep{singh2021squarm}\label{section6.1-lemma5}
Consider any two matrices $\mathbf{A} \in \mathbb{R}^{m \times n}$, $\mathbf{B} \in \mathbb{R}^{n \times n}$. Then the following holds 
 \begin{equation}
 \|\mathbf{AB}\|_F \le \|\mathbf{A}\|_F \|\mathbf{B}\|_2.
 \end{equation}
\end{Lemma}

\begin{Lemma} \citep{koloskova2019decentralized} \label{section6.1-lemma6}
For doubly stochastic matrix $\mathbf{W}$ with second largest eigenvalue $1 - \delta = |\lambda_2(\mathbf{W})| < 1$, we have
\begin{equation}
    \| \mathbf{W} - \frac{1}{n}\mathbf{1}\mathbf{1}^T \| = 1 - \delta.
\end{equation}
\end{Lemma}
\subsection{Proof of Lemmas~\ref{section3.2-lemma2-consensus error}--\ref{section3.2-lemma bound of gradient for lemma 2-3}}\label{appendix6.2-Proof of Lemmas}
For simplicity of the proof, we denote the following notations:
\begin{align}
   {{\varepsilon_1}} =~&{(1 + {\alpha _4}){{(1 - {\gamma _x}\delta )}^2}}, \nonumber\\
   {{\varepsilon_2}}=~&{(1 + \alpha _4^{ - 1})\gamma _x^2{\lambda ^2}},\nonumber\\
\varepsilon_3=~&\varepsilon_2 (1 + \alpha _1^{ - 1})((1 + \alpha_5)(1 - \omega)(1 + \alpha _6^{ - 1})+(1 + \alpha _5^{ - 1})),\nonumber\\
   \varepsilon_4=~&{\varepsilon_2}(1 + \alpha _1^{ - 1})(1 + {\alpha _5})(1 - \omega )(1 + {\alpha _6})\nonumber,\\
   {\kappa  _1} =~& (1 + \alpha _2^{ - 1})\gamma _x^2{\lambda ^2}( (1 + {\alpha _7})(1 + {\alpha _8})((1 + {\alpha _5})(1 - \omega )\nonumber\\
   &(1 + \alpha _6^{ - 1}) + (1 + \alpha _5^{ - 1})) + (1 + \alpha _7^{ - 1})),\nonumber
   \end{align}
\begin{align}
   {\kappa _2} =~& (1 + \alpha _2^{ - 1})\gamma _x^2{\lambda ^2}(1 + {\alpha _7})(1 + {\alpha _8})(1 + {\alpha _5})\nonumber\\
   &(1 - \omega )(1 + {\alpha _6}),\nonumber\\
   {\kappa _3} =~& (1 + \alpha _2^{ - 1})\gamma _x^2{\lambda ^2}(1 + {\alpha _7})(1 + \alpha _8^{ - 1}),\nonumber
\end{align}
where the parameters $\alpha_1$--$\alpha_8$ are positive constants. $\delta, \omega$ and $\lambda$ are given to us.

Based on (\ref{whh}), we rewrite 
(\ref{section3.1-algorithm1-estimete lacal momentum})--(\ref{section3.1-algorithm1-consensus}) in the following compact form to facilitate the analysis.
\begin{subequations}
\begin{align}
{\mathbf{m}_t} &= \beta {\mathbf{m}_{t - 1}} + {\mathbf{g}_t};\label{section6.2-compact form-momentum}\\
{\mathbf{x}_{t + \frac{1}{2}}} &= {\mathbf{x}_t} - \eta (\beta {\mathbf{m}_t} + {\mathbf{g}_t});\label{section6.2-compact form-local sgd update}\\
{{\hat {\mathbf{x}}}_{t + 1}} &= {{\hat{\mathbf{x}}}_t} + {\cal C}({\mathbf{x}_{t + \frac{1}{2}}} - {{\hat {\mathbf{x}}}_t});\label{section6.2-compact form-neighbors}\\
{\mathbf{x}_{t + 1}} &= {\mathbf{x}_{t + \frac{1}{2}}} + {\gamma _x}{{\hat {\mathbf{x}}}_{t + 1}}(\mathbf{W} - \mathbf{I}).\label{section6.2-compact form-consensus}
\end{align}
\end{subequations}
Now it is ready to prove \textit{Lemmas \ref{section3.2-lemma2-consensus error}--\ref{section3.2-lemma bound of gradient for lemma 2-3}}.

\textbf{1)} {\textbf{The proof of Lemma \ref{section3.2-lemma2-consensus error}}} 

We first consider the term $\| \mathbf{x}_t - \bar {\mathbf{x}}_t \|_F^2$. 

By virtue of~(\ref{section6.2-compact form-consensus}) and~\textit{Assumption~\ref{ass:graph}}, we can conclude that $\bar{\mathbf{x}}_{t}={\bar{\mathbf{x}}}_{t-\tfrac{1}{2}}$. Thus,
\begin{align}
\label{section6.2-the proof of lemma2-1}
\| \mathbf{x}_t - \bar {\mathbf{x}}_t \|_F^2 
=&\| {{\mathbf{x}_{t - \frac{1}{2}}} - {{\bar {{\mathbf{x}}}}_t} + {\gamma _x}{{\hat {\mathbf{x}}}}_t}(\mathbf{W} - \mathbf{I}) \|_F^2\nonumber \\
= &\| (\mathbf{x}_{t-\frac{1}{2}} - \bar {\mathbf{x}}_{t-\frac{1}{2}})
((1-\gamma_x)\mathbf{I} + \gamma_x \mathbf{W}) \nonumber \\
&+ \gamma_x(\hat {\mathbf{x}}_t - \mathbf{x}_{t-\frac{1}{2}})(\mathbf{W}-\mathbf{I})\|_F^2 \nonumber \\
\mathop  \le \limits^{(a)} & (1 + {\alpha _4})\|  {({\mathbf{x}_{t - \frac{1}{2}}} - {{\bar {\mathbf{x}}}_{t - \frac{1}{2}}})((1 - {\gamma _x})\mathbf{I} + {\gamma _x}\mathbf{W})} \|_F^2 \nonumber\\
 &+ (1 + \alpha _4^{ - 1})\| \gamma_x(\hat {\mathbf{x}}_t - \mathbf{x}_{t-\frac{1}{2}})(\mathbf{W}-\mathbf{I})\|_F^2\nonumber \\
 \mathop  \le \limits^{(b)} & (1 + {\alpha _4})\|  {({\mathbf{x}_{t - \frac{1}{2}}} - {{\bar {\mathbf{x}}}_{t - \frac{1}{2}}})((1 - {\gamma _x})\mathbf{I} + {\gamma _x}\mathbf{W})} \|_F^2 \nonumber\\
 &+ (1 + \alpha _4^{ - 1})\gamma _x^2\| {{{\hat {\mathbf{x}}}_t} - {{{\mathbf{x}}}_{t - \frac{1}{2}}}} \|_F^2\| {\mathbf{W} - \mathbf{I}} \|_2^2\nonumber\\
 \mathop  \le \limits^{(c)} & (1 + {\alpha _4})\|  {({\mathbf{x}_{t - \frac{1}{2}}} - {{\bar {\mathbf{x}}}_{t - \frac{1}{2}}})((1 - {\gamma _x})\mathbf{I} + {\gamma _x}\mathbf{W})} \|_F^2 \nonumber\\
 &+ (1 + \alpha _4^{ - 1})\gamma _x^2\lambda^2\| {{{\hat {\mathbf{x}}}_t} - {{{\mathbf{x}}}_{t - \frac{1}{2}}}} \|_F^2,
    \end{align}
where $(a)$, $(b)$ and $(c)$ hold due to \textit{Lemma \ref{section6.1-lemma5}} and $\lambda = \max_i \{ 1 - \lambda_i(\mathbf W) \} \Rightarrow \|\mathbf W - \mathbf I\|_2^2 \le \lambda^2$. 

Now, we analyze the two terms on the right-hand side of~(\ref{section6.2-the proof of lemma2-1}) separately. 

~~~~For the first term on the right-hand side of (\ref{section6.2-the proof of lemma2-1}), we have
\begin{align}
\label{section6.2-the proof of lemma2-2}
\|  {({\mathbf{x}_{t - \frac{1}{2}}}} &-{ {{\bar {\mathbf{x}}}_{t - \frac{1}{2}}})}{[(1 - {\gamma _x})\mathbf{I} + {\gamma _x}\mathbf{W}]} \|_F\nonumber\\
\le & (1 - {\gamma _x}){\| {{\mathbf{x}_{t - \frac{1}{2}}} - {{\bar {{\mathbf{x}}}}_{t - \frac{1}{2}}}} \|_F} \nonumber+ {\gamma _x}{\| {({\mathbf{x}_{t - \frac{1}{2}}} - {{\bar {{\mathbf{x}}}}_{t - \frac{1}{2}}}){\bf{W}}} \|_F}\nonumber\\
=&(1 - {\gamma _x}){\| {{{{\mathbf{x}}}_{t - \frac{1}{2}}} - {{\bar {{\mathbf{x}}}}_{t - \frac{1}{2}}}} \|_F} \nonumber\\
&+ {\gamma _x}{\| {({{{\mathbf{x}}}_{t - \frac{1}{2}}} - {{\bar {{\mathbf{x}}}}_{t - \frac{1}{2}}})({\bf{W}} - \frac{{{{\mathbf{11}}^T}}}{n})} \|_F}\nonumber\\
\le & (1 - {\gamma _x}){\| {{{{\mathbf{x}}}_{t - \frac{1}{2}}} - {{\bar {{\mathbf{x}}}}_{t - \frac{1}{2}}}} \|_F}+ {\gamma _x}(1-\delta) {\| {{{{\mathbf{x}}}_{t - \frac{1}{2}}} - {{\bar {{\mathbf{x}}}}_{t - \frac{1}{2}}}} \|_F}\nonumber\\
= & (1 - {\gamma _x}\delta ){\| {{{{\mathbf{x}}}_{t - \frac{1}{2}}} - {{\bar {{\mathbf{x}}}}_{t - \frac{1}{2}}}} \|_F},
\end{align}
 where the first equality holds due to $({{{{\mathbf{x}}}}_{t - \frac{1}{2}}}-{{\bar {{\mathbf{x}}}}_{t - \frac{1}{2}}})\frac{{{{\mathbf{11}}^T}}}{n} = \mathbf{0}$ and the last inequality holds based on \textit{Lemma \ref{section6.1-lemma6}}.

~~~~Now, we show the upper bound of $\| {{\mathbf{x}_{t - \frac{1}{2}}} - {{\bar {{\mathbf{x}}}}_{t - \frac{1}{2}}}} \|_F^2$. 
Taking the average of both sides of (\ref{section6.2-compact form-local sgd update}), we have
 \begin{equation}\label{section6.2-the proof of lemma2-3}
{\mathbf{\bar x}_{t-\frac{1}{2}}} = {\mathbf{\bar x}_{t-1}} - \eta (\beta {\mathbf{m}_{t-1}} + {\mathbf{g}_{t-1}})\frac{{{{\mathbf{11}}^T}}}{n},
 \end{equation} 
From~(\ref{section6.2-compact form-local sgd update}) and (\ref{section6.2-the proof of lemma2-3}), we have
\begin{align}\label{section6.2-the proof of lemma2-4}
&\mathbb{E}[\|{\mathbf{x}_{t - \frac{1}{2}}} - {{\bar {{\mathbf{x}}}}_{t - \frac{1}{2}}}\|_F^2]\nonumber\\
&\le \mathbb{E}[\|{{{\mathbf{x}}}_{t - 1}} -{{ { \bar{\mathbf{x}}}}_{t - 1}} - \eta (\beta {\mathbf{m}_{t - 1}} + {\mathbf{g}_{t - 1}})(\frac{{{{\mathbf{11}}^T}}}{n}-\mathbf{I})\|_F^2].   
\end{align}
Similarly, from~(\ref{section6.2-compact form-local sgd update}), we have
 \begin{equation}\label{section6.2-the proof of lemma2-5}
\mathbb{E}[\| {{{\hat {\mathbf{x}}}_t} - {{{\mathbf{x}}}_{t - \frac{1}{2}}}}\| _F^2]\le \mathbb{E}[\| {{{\hat {\mathbf{{x}}}}_t} - {\mathbf{x}_{t - 1}} + \eta (\beta {\mathbf{m}_{t - 1}} + {\mathbf{g}_{t - 1}})} \|_F^2].
 \end{equation}
 By substituting~(\ref{section6.2-the proof of lemma2-2}), (\ref{section6.2-the proof of lemma2-3}), (\ref{section6.2-the proof of lemma2-5}) into~(\ref{section6.2-the proof of lemma2-1}) and taking the expectation of both sides yields that
\begin{align} 
\label{section6.2-the proof of lemma2-6}
\mathbb{E}[\| \mathbf{x}_t &- \bar {\mathbf{x}}_t \|_F^2] \nonumber\\
 \le&  {{\varepsilon_1}}\mathbb{E}[\| {{\mathbf{x}_{t - \frac{1}{2}}} - {{\bar {\mathbf{x}}}_{t - \frac{1}{2}}}} \|_F^2] + {{\varepsilon_2}}\mathbb{E}[\| {{{\hat {\mathbf{x}}}_t} - {{ \mathbf{x}}_{t - \frac{1}{2}}}} \|_F^2] \nonumber\\
  \le& {\varepsilon_1}\mathbb{E}[\| {{{{\mathbf{x}}}_{t - 1}} -{\bar  {\mathbf{x}}_{t-1}} - \eta (\beta {\mathbf{m}_{t - 1}} + {\mathbf{g}_{t - 1}})(\frac{{\mathbf{1}{\mathbf{1}^T}}}{n} - \mathbf{I})} \|_F^2] \nonumber\\
  &+ {\varepsilon_2}\mathbb{E}[\| {{{\hat {\mathbf{x}}}_t} - {{ \mathbf{x}}_{t -1}+\eta (\beta {\mathbf{m}_{t - 1}} + {\mathbf{g}_{t - 1}})}}\|_F^2],\nonumber\\
 \le & {\varepsilon _1}(1 + \alpha _1^{ - 1})\mathbb{E}[\| {{{ {\mathbf{x}}}_{t - 1}} - {\bar{\mathbf{x}}_{t - 1}}} \|_F^2] \nonumber\\
 &+ {\varepsilon _1}(1 + {\alpha _1})\mathbb{E}[\| {\eta (\beta {{\mathbf{m}}_{t - 1}} + {{\mathbf{g}}_{t - 1}})(\frac{{{{\mathbf{11}}^T}}}{n} - \mathbf{I})}\|_F^2]\nonumber\\
 &+ {\varepsilon _1}(1 + \alpha _1^{ - 1})\mathbb{E}[\| {{{\hat {\mathbf{x}}}_t} - {{\mathbf{x}}_{t - 1}}}\|_F^2] \nonumber\\
 &+ {\varepsilon _2}(1 + {\alpha _1})\mathbb{E}[\| {\eta (\beta {{\mathbf{m}}_{t - 1}} + {{\mathbf{g}}_{t - 1}})}\|_F^2]\nonumber\\
\le & {\varepsilon _1}(1 + \alpha _1^{ - 1})\mathbb{E}[\| {{{{\mathbf{x}}}_{t - 1}} - {\bar {\mathbf{x}}_{t - 1}}} \|_F^2]\nonumber\\
&+ {\varepsilon _2}(1 + \alpha _1^{ - 1})\mathbb{E}[\| {{{\hat {\mathbf{x}}}_t} - {{\mathbf{x}}_{t - 1}}} \|_F^2]\nonumber\\
& + (1 + {\alpha _1})({\varepsilon _1} + {\varepsilon _2})\mathbb{E}[\| {\eta (\beta {{\mathbf{m}}_{t - 1}} + {{\mathbf{g}}_{t - 1}})} \|_F^2],
   \end{align}
where the last equality holds due to ${\| {\frac{{\mathbf{1}{\mathbf{1}^T}}}{n} - \mathbf{I}} \|_2} = 1$.

An important observation is that 
\begin{align}
\label{section6.2-the proof of lemma2-7}
 \mathbb{E}[\| {{{\hat {\mathbf{x}}}_t}} &- {{\mathbf{x}_{t - 1}}} \|_F^2] \nonumber\\
 \mathop =\limits^{(a)}&\mathbb{E}{[\| {{\mathbf{x}}_{t - 1}}-({{\hat {\mathbf{x}}}_{t - 1}}+{\cal C}({{\mathbf{x}}_{t - \frac{1}{2}}} - {{\hat {\mathbf{x}}}_{t - 1}}))\|_F^2]}\nonumber\\
=&\mathbb{E}{[\| {{{\mathbf{x}}_{t - \frac{1}{2}}} - {{\hat {\mathbf{x}}}_{t - 1}} - {\cal C}({{\mathbf{x}}_{t - \frac{1}{2}}} - {{\hat {\mathbf{x}}}_{t - 1}}) - {{\mathbf{x}}_{t - \frac{1}{2}}} + {{\mathbf{x}}_{t - 1}}} \|_F^2}] \nonumber\\
\mathop  \le \limits^{(b)} &{\rm{(1 + }}{\alpha _5}{\rm{)(1}} - \omega {\rm{)\mathbb{E}}}{[\| {{\mathbf{x}_{t - \frac{1}{2}}} - {{\hat {\mathbf{x}}}_{t - 1}}} \|_F^2}]\nonumber\\
&+ (1 + \alpha _5^{ - 1})\mathbb{E}{[\| {{\mathbf{x}_{t - 1}} - {\mathbf{x}_{t - \frac{1}{2}}}} \|_F^2}]\nonumber \\
=&{\rm{(1 + }}{\alpha _5}{\rm{)(1}} - \omega {\rm{)\mathbb{E}}}[\| {{\mathbf{x}_{t - \frac{1}{2}}} -{{\mathbf{x}}}_{t - 1}}+{{\mathbf{x}}}_{t - 1}+{{\hat {\mathbf{x}}}_{t - 1}} \|_F^2]\nonumber\\
&+ (1 + \alpha _5^{ - 1})\mathbb{E}{[\| {{\mathbf{x}_{t - 1}} - {\mathbf{x}_{t - \frac{1}{2}}}} \|_F^2}] \nonumber\\
\mathop  \le \limits^{(c)} &{\rm{(1 + }}{\alpha _5}{\rm{)(1}} - \omega {\rm{)[}}(1 + \alpha _6^{ - 1}){\rm{\mathbb{E}}}{[\| {{\mathbf{x}_{t - \frac{1}{2}}} - {\mathbf{x}_{t - 1}}} \|_F^2}]  \nonumber\\
 &+ {\rm{(1 + }}{\alpha _6}{\rm{)\mathbb{E}}}{[\| {{\mathbf{x}_{t - 1}} - {{\hat {\mathbf{x}}}_{t - 1}}} \|_F^2}] \nonumber\\
 &+ (1 + \alpha _5^{ - 1})\mathbb{E}{[\| {{\mathbf{x}_{t - 1}} - {{\mathbf{x}}_{t - \frac{1}{2}}}} \|_F^2}] \nonumber\\ 
 =& ((1 + {\alpha _5})(1 - \omega )(1 + \alpha _6^{ - 1}) + (1 + \alpha _5^{ - 1}))\nonumber\\ 
 &\mathbb{E}{[\| {{{\mathbf{x}}_{t - \frac{1}{2}}} - {{\mathbf{x}}_{t - 1}}} \|_F^2}]+ (1 + {\alpha _5})(1 - \omega ){\rm{(1 + }}{\alpha _6})\nonumber\\
 &\mathbb{E}{[\| {{{\mathbf{x}}_{t - 1}} - {{\hat {\mathbf{x}}}_{t - 1}}} \|_F^2]},
    \end{align}
where $(a)$ holds due to
(\ref{section6.2-compact form-neighbors}); $(b)$ and $(c)$ hold due to the \textit{Definition~\ref{section2.4-definition-1}} and \textit{Lemma~\ref{section6.1-lemma5}}, respectively.

From~(\ref{section6.2-compact form-local sgd update}), we have
\begin{equation}\label{section6.2-lemma2-1111111}
    \mathbb{E}{[\| {{{\mathbf{x}}_{t - \frac{1}{2}}} - {{\mathbf{x}}_{t - 1}}} \|_F^2}]=\mathbb{E}{[\|\eta (\beta {\mathbf{m}_{t-1}} + {\mathbf{g}_{t-1}})\|_F^2}].
\end{equation}
Consequently, substituting~(\ref{section6.2-the proof of lemma2-7}) into~(\ref{section6.2-the proof of lemma2-6}), we have
\begin{align}\label{section6.2-the proof of lemma2-8}
\mathbb{E}[\| \mathbf{x}_t &- \bar {\mathbf{x}}_t \|_F^2]\nonumber 
\le {\varepsilon _1}(1 + \alpha _1^{ - 1})\mathbb{E}[\| {{{{\mathbf{x}}}_{t - 1}} - {\bar {\mathbf{x}}_{t - 1}}} \|_F^2]\nonumber\\
&+{\varepsilon _3}{{\mathbb{E}}{[\| {{{\mathbf{x}}_{t - \frac{1}{2}}} - {{\mathbf{x}}_{t - 1}}} \|_F^2]}}+{{{{\varepsilon _4}\mathbb{E}}}{[\| {{{\mathbf{x}}_{t - 1}} - {{\hat {\mathbf{x}}}_{t - 1}}} \|_F^2]}}\nonumber\\
&+ (1 + {\alpha _1})({\varepsilon _1} + {\varepsilon _2})\mathbb{E}[\| {\eta (\beta {{\mathbf{m}}_{t - 1}} + {{\mathbf{g}}_{t - 1}})} \|_F^2].
   \end{align}
From~(\ref{section6.2-lemma2-1111111}),~(\ref{section6.2-the proof of lemma2-8}) can be rearranged as~(\ref{section3.2-lemma2-consensus error-1}).

\textbf{2)} {\textbf{The proof of Lemma 3}}

From~(\ref{section6.2-compact form-neighbors})--(\ref{section6.2-compact form-consensus}), we have
\begin{align}
\label{section6.2-the proof of lemma3-0}
 \mathbb{E}[\| {{{\mathbf{x}}_t}} &- {{{\hat {\mathbf{x}}}_t}} \|_F^2]\nonumber=\mathbb{E}[\| {{{\mathbf{x}}_t} - ({{\hat {\mathbf{x}}}_{t - 1}} + {\cal C}({{\mathbf{x}}_{t - \frac{1}{2}}} - {{\hat {\mathbf{x}}}_{t - 1}}))} \|_F^2] \nonumber\\
  =& \mathbb{E}[\| {{{\mathbf{x}}_{t - \frac{1}{2}}} - {{\hat {\mathbf{x}}}_{t - 1}} - {\cal C}({{\mathbf{x}}_{t - \frac{1}{2}}} - {{\hat {\mathbf{x}}}_{t - 1}})+{\mathbf{x}}_t} -{{\mathbf{x}}_{t - \frac{1}{2}}} \|_F^2]\nonumber\\
  \le& (1 + {\alpha _2})(1 - \omega )\mathbb{E}[\| {{{\mathbf{x}}_{t - \frac{1}{2}}} - {{\hat {\mathbf{x}}}_{t - 1}}}\|_F^2]\nonumber\\
  &+ (1 + \alpha _2^{ - 1})\mathbb{E}[\| {{{\mathbf{x}}_t} - {{\mathbf{x}}_{t - \frac{1}{2}}}} \|_F^2]. 
 \end{align}
~~~~Now, we provide the upper bound of $\mathbb{E}[\| {{{\mathbf{x}}_t} - {{\mathbf{x}}_{t - \frac{1}{2}}}} \|_F^2]$. From~(\ref{section6.2-compact form-consensus}), we have
\begin{align}
\label{section6.2-the proof of lemma3-1}
 \mathbb{E}[\| {{{\mathbf{x}}_t}} &{- {{\mathbf{x}}_{t - \frac{1}{2}}}}\|_F^2]\nonumber= \gamma _x^2\mathbb{E} [\| {{{\hat {\mathbf{x}}}_t}(\mathbf{W} - \mathbf{I})} \|_F^2] \nonumber\\ 
 \mathop  = \limits^{(a)}  &\gamma _x^2\mathbb{E}[\| {({{\hat {\mathbf{x}}}_t} - {{\bar {\mathbf{x}}}_{t - \frac{1}{2}}})(\mathbf{W} - \mathbf{I})} \|_F^2] \nonumber\\
 \mathop  \le \limits^{(b)} &\gamma _x^2{\lambda ^2}\mathbb{E}[\| {{{\hat {\mathbf{x}}}_t} - {{\bar {\mathbf{x}}}_{t - \frac{1}{2}}}} \|_F^2]\nonumber \\
 =& \gamma _x^2{\lambda ^2}\mathbb{E}[\| {{{\hat {\mathbf{x}}}_t} - {\mathbf{\bar x}_{t-1}} + \eta (\beta {\mathbf{m}_{t-1}} + {\mathbf{g}_{t-1}})\frac{{{{\mathbf{11}}^T}}}{n}} \|_F^2]  \nonumber\\
 \le &\gamma _x^2{\lambda ^2}(1 + {\alpha _7}) \mathbb{E}[\| {{{\hat {\mathbf{x}}}_t} - {\mathbf{\bar x}_{t-1}}} \|_F^2]+\gamma _x^2{\lambda ^2}(1 + \alpha _7^{ - 1})\nonumber\\
 &\mathbb{E}[\| { \eta (\beta {\mathbf{m}_{t-1}} + {\mathbf{g}_{t-1}})\frac{{{{\mathbf{11}}^T}}}{n}} \|_F^2] \nonumber\\
=&\gamma _x^2{\lambda ^2}(1 + {\alpha _7}) \mathbb{E}[\| {{{\hat {\mathbf{x}}}_t} - {\mathbf{\bar x}_{t-1}}} \|_F^2]+\gamma _x^2{\lambda ^2}(1 + \alpha _7^{ - 1})\nonumber\\
 &\mathbb{E}[\| { \eta (\beta {\mathbf{m}_{t-1}} + {\mathbf{g}_{t-1}}}) \|_F^2]\nonumber\\
\mathop  \le \limits^{(c)} & \gamma _x^2{\lambda ^2}(1 + {\alpha _7})({(1 + {\alpha _8})\mathbb{E}[\| {{{\hat {\mathbf{x}}}_t} - {{{\mathbf{x}}}_{t - 1}}} \|_F^2})\nonumber\\
 &+ (1 + \alpha _8^{ - 1})\mathbb{E}[\| {{{\mathbf{x}}_{t - 1}} - {{\bar {\mathbf{x}}}_{t - 1}}} \|_F^2] \nonumber\\
 &+ \gamma _x^2{\lambda ^2}(1 + \alpha _7^{ - 1})\mathbb{E}[\| { \eta (\beta {\mathbf{m}_{t-1}} + {\mathbf{g}_{t-1}})} \|_F^2] \nonumber\\
 =&\gamma _x^2{\lambda ^2}{(1 + {\alpha _7})(1 + {\alpha _8})\mathbb{E}[\| {{{\hat {\mathbf{x}}}_t} - {{{\mathbf{x}}}_{t - 1}}} \|_F^2}]\nonumber\\
& +\gamma _x^2{\lambda ^2}(1 + {\alpha _7})(1 + \alpha _8^{ - 1})\mathbb{E}[\| {{{\mathbf{x}}_{t - 1}} - {{\bar {\mathbf{x}}}_{t - 1}}} \|_F^2]\nonumber\\
 &+ \gamma _x^2{\lambda ^2}(1 + \alpha _7^{ - 1})\mathbb{E}[\| { \eta (\beta {\mathbf{m}_{t-1}} + {\mathbf{g}_{t-1}})} \|_F^2],
 \end{align}
where $(a)$ holds due to ${{\bar {\mathbf{x}}}_{t - \frac{1}{2}}}(\mathbf{W} - \mathbf{I}) = \mathbf{0}$. The inequality $(b)$ follows $\lambda = \max_i \{ 1 - \lambda_i(\mathbf W) \} \Rightarrow \|\mathbf W - \mathbf I\|_2^2 \le \lambda^2$ and the validity of $(c)$ follows from \textit{Lemma~\ref{section6.1-lemma5}}.

By combining~(\ref{section6.2-the proof of lemma3-1}) with~(\ref{section6.2-the proof of lemma2-7}), we have 
\begin{align}\label{section6.2-the proof of lemma3-2}
\mathbb{E}&[\| {{{\mathbf{x}}_t} - {{\mathbf{x}}_{t - \frac{1}{2}}}} \|_F^2]\nonumber\\
\le& {\gamma _x^2}{\lambda ^2}{(1 + {\alpha _7})(1 + {\alpha _8})((1 + {\alpha _5})(1 - \omega )(1 + \alpha _6^{ - 1}) }\nonumber\\
&+ {(1 + \alpha _5^{ - 1}){\rm{)\mathbb{E}}}{[\| {{{\mathbf{x}}_{t - \frac{1}{2}}} - {{\mathbf{x}}_{t - 1}}} \|_F^2]}}\nonumber\\
&+ {\gamma _x^2{\lambda ^2}{(1 + {\alpha _7})(1 + {\alpha _8})}} (1 + {\alpha _5})(1 - \omega )\nonumber\\
&~~~{\rm{(1 + }}{\alpha _6}{\rm{)\mathbb{E}}}{[\| {{{\mathbf{x}}_{t - 1}} - {{\hat {\mathbf{x}}}_{t - 1}}} \|_F^2]}\nonumber\\
& +\gamma _x^2{\lambda ^2}(1 + {\alpha _7})(1 + \alpha _8^{ - 1})\mathbb{E}[\| {{{\mathbf{x}}_{t - 1}} - {{\bar {\mathbf{x}}}_{t - 1}}} \|_F^2]\nonumber\\
 &+ \gamma _x^2{\lambda ^2}(1 + \alpha _7^{ - 1})\mathbb{E}[\| { \eta (\beta {\mathbf{m}_{t-1}} + {\mathbf{g}_{t-1}})} \|_F^2].
\end{align}
~~~Then, the upper bound of $\mathbb{E}[\| {{{\mathbf{x}}_{t - \frac{1}{2}}} - {{\hat {\mathbf{x}}}_{t - 1}}}\|_F^2]$ is analyzed. Based on \textit{Lemma~\ref{section6.1-lemma5}}, it readily follows that
\begin{align}\label{lemma 3-i-4}
 \mathbb{E}&[\| {{{\mathbf{x}}_{t - \frac{1}{2}}} - {{\hat {\mathbf{x}}}_{t - 1}}}\|_F^2]\nonumber\\ 
 =&\mathbb{E}[\| {{{\mathbf{x}}_{t - 1}} - \eta (\beta \mathbf{m}_{t - 1}} + {\mathbf{g}_{t - 1}}) - {{\hat {\mathbf{x}}}_{t - 1}}\|_F^2] \nonumber\\ 
\le& (1 + {\alpha _3})\mathbb{E}[\| {{{\mathbf{x}}_{t - 1}} - {{\hat {\mathbf{x}}}_{t - 1}}} \|_F^2]\nonumber\\ 
  &+ (1 + \alpha _3^{ - 1}) \mathbb{E}[\| {\eta(\beta {\mathbf{m}_{t - 1}} + {\mathbf{g}_{t - 1}}})\|_F^2]. 
 \end{align}
By combining~(\ref{section6.2-the proof of lemma3-1}) with~(\ref{lemma 3-i-4}), (\ref{section6.2-the proof of lemma3-0}) can be rearranged as~(\ref{section3.2-lemma2-compressed error-1}).

\textbf{3)} {\textbf{The proof of Lemma 4}} 

To begin with, we denote following notations:
\begin{align}
\label{light}
{\ell_1}=~&{\varepsilon_1(1+\alpha_1^{-1})+\kappa_3},\nonumber \\
{\ell_2}=~&{\varepsilon_4+(1+\alpha_2)(1-\omega)(1+\alpha_3)+\kappa_2},\nonumber\\
 {{\ell _3}}=~&{\varepsilon _5}+(1 + {\alpha _1})({\varepsilon _1} + {\varepsilon _2})+(1 + {\alpha _2})(1 - \omega )\nonumber\\
 &(1 + {\alpha _3^{-1}})+{\kappa _1},\nonumber\\
 {\kappa _4}=~&\frac{{4n{d^2}\sigma _2^2}(L_{{f_1}}^2\gamma _g^2\sigma _2^2 + 2\gamma _1 + 2{\vartheta _1})}{{\gamma _g^2}}. 
\end{align}
We define an auxiliary sequence
${\tilde x_{i,t}}$ for each agent $i$ as follows:
\begin{equation}\label{section6.2-lemma4-auxiliary sequence}
{\tilde x_{i,t}} = {x_{i,t}} - \frac{{\eta {\beta ^{\rm{2}}}}}{{1 - \beta }}{m_{i,t - 1}},
\end{equation}
Let ${\bar x}_t^a = \frac{1}{n}\sum\limits_{i = 1}^n {{{\tilde x}_{i,t}}}$. 
From~(\ref{section6.2-lemma4-auxiliary sequence}), we have
\begin{equation}\label{section6.2-lemma4-the average of vauxiliary sequence}
{\bar x}_t^a = {{\bar x}_t} - \frac{{\eta {\beta ^{\rm{2}}}}}{{1 - \beta }}\frac{1}{n}\sum\limits_{i = 1}^n {{m_{i,t - 1}}}.
\end{equation}
From~(\ref{section6.2-lemma4-auxiliary sequence})--(\ref{section6.2-lemma4-the average of vauxiliary sequence}), we have
\begin{align}
\label{section6.2-lemma4-recurrence relation of the average of auxiliary sequence}
{\bar x}_{t+1}^a &= {{\bar x}_{t + 1}} - \frac{{\eta {\beta ^{\rm{2}}}}}{{1 - \beta }}\frac{1}{n}\sum\limits_{i = 1}^n {{m_{i,t}}} \nonumber \\
 &= {{\bar x}_t} - \frac{\eta }{n}\sum\limits_{i = 1}^n {(\beta {m_{i,t}} + {g_{i,t}})}  - \frac{{\eta {\beta ^{\rm{2}}}}}{{1 - \beta }}\frac{1}{n}\sum\limits_{i = 1}^n {{m_{i,t}}}  \nonumber\\
 &= {{\bar x}_t} - \frac{\eta }{n}\sum\limits_{i = 1}^n {{g_{i,t}}}  - \frac{1}{n}\left(\eta \beta  + \frac{{\eta {\beta ^{\rm{2}}}}}{{1 - \beta }}\right)\sum\limits_{i = 1}^n {{m_{i,t}}} \nonumber\\
 & ={{\bar x}_t} - \frac{\eta }{n}\sum\limits_{i = 1}^n {{g_{i,t}}}  - \frac{1}{n}\frac{{\eta \beta }}{{1 - \beta }}\sum\limits_{i = 1}^n {(\beta {m_{i,t - 1}} + {g_{i,t}})} \nonumber\\
& = {{\bar x}_t} - \left(\eta  + \frac{{\eta \beta }}{{1 - \beta }}\right)\frac{1}{n}\sum\limits_{i = 1}^n {{g_{i,t}}}  - \frac{{\eta {\beta ^2}}}{{1 - \beta }}\frac{1}{n}\sum\limits_{i = 1}^n {{m_{i,t-1}}}\nonumber\\
 &={\bar x}_t^a - \frac{\eta }{{1 - \beta }}\frac{1}{n}\sum\limits_{i = 1}^n {{g_{i,t}}},
  \end{align}
where the second and fourth equalities hold due to~(\ref{section3.1-algorithm1-local SGD Update}) and (\ref{section3.1-algorithm1-estimete lacal momentum}), respectively.

Then, from (\ref{section3.2-Lyapunov candidate function}), we have 
\begin{equation}
\label{section6.2-lemma4-Xt express}
\begin{aligned}
{\chi _t} =& \mathbb{E}[\| {{\mathbf{x}_t} - {{\bar {\mathbf{x}}}_t}} \|_F^2] + \mathbb{E}[\| {{\mathbf{x}_t} - {{\hat {\mathbf{x}}}_t}} \|_F^2]\\
\le &({\varepsilon _1}(1 + \alpha _1^{ - 1})+{\kappa _3})\mathbb{E}[\| {{{ {\mathbf{x}}}_{t - 1}} - {\bar{\mathbf{x}}_{t - 1}}} \|_F^2]\\
&+({\varepsilon _3}+(1 + {\alpha _1})({\varepsilon _1} + {\varepsilon _2})+(1 + {\alpha _2})(1 - \omega )(1 + {\alpha _3^{ - 1}})\\
&+{\kappa _1})\mathbb{E}[\| {\eta (\beta {{\mathbf{m}}_{t - 1}} + {{\mathbf{g}}_{t - 1}})} \|_F^2]+({\varepsilon _4}+(1 + {\alpha _2})(1- \omega )\\
&~~~(1 + {\alpha _3})+{\kappa _2})\mathbb{E}[\| {{{\mathbf{x}}_{t - 1}} - {{\hat {\mathbf{x}}}_{t - 1}}} \|_F^2]\\
\le& \max\{{{\ell_1}, {\ell_2}}\}{\chi _{t-1}}+{\ell _3}\mathbb{E}[\| {\eta (\beta {{\mathbf{m}}_{t - 1}} + {{\mathbf{g}}_{t - 1}})} \|_F^2],
\end{aligned}
\end{equation}
where the first inequality holds due to \textit{Lemmas \ref{section3.2-lemma2-consensus error}} and \textit{Lemma \ref{section3.2-lemma2-compressed error}}.

Now, we analyze the bound of $\mathbb{E}[\| {\eta (\beta {{\mathbf{m}}_{t - 1}} + {{\mathbf{g}}_{t - 1}})} \|_F^2]$. 

For all $t \ge 0$, we first consider the term $\mathbb{E}[\| {{\mathbf{g}_t}} \|_F^2]$. From (\ref{section2.3-one point in this paper}), we have
\begin{align*} 
\mathbb{E}&_{{{\cal L}_t}}[\| {{\mathbf{g}_t}} \|_F^2]\nonumber
= \sum\limits_{i = 1}^n {\mathbb{E}_{{{\cal L}_t}}{{[ \| {{g_{i,t}}} \|}^2]}}\nonumber \\
 =& \sum\limits_{i = 1}^n {\mathbb{E}_{{{\cal L}_t}}{{ [\| {\frac{d}{{{\gamma _g}}}({{f_i}({x_{i,t}} + {\gamma _g}{u_{i,t}},{\xi _{i,t}}) + {\varphi _{i,t}}}){u_{i,t}}\|^2]}}}}\nonumber
 \end{align*} 
 \begin{align} 
\label{section6.2-lemma4-bounded of the expected norm squared}
 =& \frac{{{d^2}}}{{\gamma _g^2}}\sum\limits_{i = 1}^n {\mathbb{E}_{{{\cal L}_t}} [| {{f_i}({x_{i,t}} + {\gamma _g}{u_{i,t}},{\xi _{i,t}}) - {f_i}({{\bar x}_t} + {\gamma _g}{u_{i,t}},{\xi _{i,t}})}} \nonumber\\
 &+ {f_i}({{\bar x}_t} + {\gamma _g}{u_{i,t}},{\xi _{i,t}}) - {f_i}({{\bar x}_t},{\xi _{i,t}}) \nonumber\\
 &+ {f_i}({{\bar x}_t},{\xi _{i,t}}){ { - {f_i}({x_{i,t}},{\xi _{i,t}}) + {f_i}({x_{i,t}},{\xi _{i,t}})+ {\varphi _{i,t}}|^2\|{u_{i,t}}}  \|^2]}\nonumber\\
 \le& \frac{{{d^2}}}{{\gamma _g^2}}\sum\limits_{i = 1}^n {\mathbb{E}_{{{\cal L}_t}} [4(| {{f_i}({x_{i,t}} + {\gamma _g}{u_{i,t}},{\xi _{i,t}}) - {f_i}({{\bar x}_t} + {\gamma _g}{u_{i,t}},{\xi _{i,t}})|^2}} \nonumber\\
 &+ |{f_i}({{\bar x}_t} + {\gamma _g}{u_{i,t}},{\xi _{i,t}}) - {f_i}({{\bar x}_t},{\xi _{i,t}})|^2 \nonumber\\
 &+ |{f_i}({{\bar x}_t},{\xi _{i,t}}){ { - {f_i}({x_{i,t}},{\xi _{i,t}})|^2  +| {f_i}({x_{i,t}},{\xi _{i,t}}) + {\varphi _{i,t}}|^2 )\|{u_{i,t}}}  \|^2]}\nonumber\\
\le& \frac{{4{d^2}}}{{\gamma _g^2}}\sum\limits_{i = 1}^n {\mathbb{E}_{{{\cal L}_t}}[{{ \| {{u_{i,t}}}  \|}^2}(2L_{{f_1}}^2{\|x_{i,t}} - {{{\bar x}_t}\|{^2} + L_{{f_1}}^2\gamma _g^2\|{u_{i,t}}\|{^2}}}\nonumber\\
 &+2|{f_i}({x_{i,t}},{\xi _{i,t}})|{^2} +2|{\varphi _{i,t}}|^2)]\nonumber\\
\le & \frac{{8{d^2}\sigma _2^2 L_{{f_1}}^2}}{{\gamma _g^2}}\sum\limits_{i = 1}^n \|{{x_{i,t}} - {{\bar x}_t}\|{^2}}  + \kappa _4,
\end{align}
where the second and last inequalities hold due to \textit{Assumptions \ref{section 2.1-ass:the local function},~\ref{section2.3-ass:noise}} and~\textit{\ref{section2.3-ass:perturbation vector}}. 

Then, for every $t\ge0$, taking expectation on both sides of (\ref{section6.2-lemma4-bounded of the expected norm squared}) yields
\begin{align} 
\label{0990}
\mathbb{E} [\| {{\mathbf{g}_t}} \|_F^2]
\le & \frac{{8{d^2}\sigma _2^2 L_{{f_1}}^2}}{{\gamma _g^2}}\sum\limits_{i = 1}^n \mathbb{E}[\|{{x_{i,t}} - {{\bar x}_t}\|{^2}]}  + \kappa _4\nonumber\\
\le & \frac{{8{d^2}\sigma _2^2 L_{{f_1}}^2}}{{\gamma _g^2}}\chi_t + \kappa _4,
\end{align}
where the second inequality holds due to (\ref{section3.2-Lyapunov candidate function}). 

Let ${{\bar g}_t}={\frac{1}{n}\sum\limits_{i = 1}^n {{g_{i,t}}} }$, based on Jensen's inequality, 
we have
\begin{equation}\label{section6.2-lemma4-gt average norm square}
    {\| {{{\bar g}_t}} \|^2} = {\| {\frac{1}{n}\sum\limits_{i = 1}^n {{g_{i,t}}} } \|^2} \le \frac{1}{n}{\sum\limits_{i = 1}^n {\| {{g_{i,t}}} \|} ^2}.
\end{equation}
From (\ref{section6.2-lemma4-bounded of the expected norm squared}), we have
\begin{align}
\label{section6.2-lemma4-gt average norm square and take expectation}
\mathbb{E}_{{{\cal L}_t}}{[\| {{{\bar g}_t}} \|^2]}&\le\frac{1}{n}{\sum\limits_{i = 1}^n \mathbb{E}_{{{\cal L}_t}}[{\| {{g_{i,t}}} \|} ^2}]\nonumber\\
&\le \frac{{8{d^2}\sigma _2^2 L_{{f_1}}^2}}{{n\gamma _g^2}}\sum\limits_{i = 1}^n \|{{x_{i,t}} - {{\bar x}_t}\|{^2}} + \frac{{{\kappa _4}}}{n}.
\end{align}

Then, we consider the term $\mathbb{E}[\| {{\mathbf{m}_{t-1}}} \|_F^2]$.

From the~(\ref{section6.2-compact form-momentum}) and $\mathbf{m}_{-1}=\mathbf{0}$, we can derive that
\begin{equation}\label{section6.2-lemma4-mt-1 expression}
\mathbf{m}_{t-1}
= \beta^{t}\mathbf{m}_{-1}+\sum_{k=0}^{t-1}\beta^{t-1-k}\mathbf{g}_{k}=\sum_{k=0}^{t-1}\beta^{t-1-k}\mathbf{g}_{k}.
\end{equation}
Denote $\theta_{t-1}=\sum_{k=0}^{t-1}\beta^{t-1-k}$. Due to $\beta \in [0,1)$, it follows that 
\begin{align}\label{111111}
  \theta_{t-1}=\frac{1-\beta^{t}}{1-\beta}\le\frac{1}{1-\beta}.  
\end{align}
Let ${\alpha _k} = \frac{{{\beta ^{t - 1 - k}}}}{{{\theta _{t - 1}}}}$, we have $\sum\limits_{k = 0}^{t - 1} {{\alpha _k}}= 1$. From (\ref{section6.2-lemma4-mt-1 expression}), we have 
\begin{align*}
\mathbb{E}_{{{\cal L}_{t-1}}}[\| {{{\bf{m}}_{t - 1}}} \|_F^2]  &= \mathbb{E}_{{{\cal L}_{t-1}}}\left[\left\| {\sum\limits_{k = 0}^{t - 1} {{\beta ^{t - 1 - k}}} {{\bf{g}}_k}} \right\|_F^2\right]\nonumber\\
&= \mathbb{E}_{{{\cal L}_{t-1}}}\left[\left\| {{\theta _{t - 1}}\sum\limits_{k = 0}^{t - 1} {\frac{{{\beta ^{t - 1 - k}}}}{{{\theta _{t - 1}}}}} {{\bf{g}}_k}} \right\|_F^2\right]\nonumber
\end{align*}
 \begin{align}
 \label{section6.2-the proof of m_t-1 norm expectation}
& = \theta _{t - 1}^2 \mathbb{E}_{{{\cal L}_{t-1}}}\left[\left\| {\sum\limits_{k = 0}^{t - 1} {\alpha _k} {{\bf{g}}_k}}  \right\|_F^2\right]\nonumber\\
&\le \theta_{t-1}^2 \sum_{k=0}^{t-1} \alpha_k \mathbb{E}_{{{\cal L}_{t-1}}}[\|{{\bf{g}}_k}\|_F^2]\nonumber\\
 &= {\theta _{t - 1}}\sum\limits_{k = 0}^{t - 1} {{\beta ^{t - 1 - k}}} \mathbb{E}_{{{\cal L}_{t-1}}}[\| {{{\bf{g}}_k}} \|_F^2].
\end{align}
where the inequality holds due to Jensen's inequality. 

Taking expectation on both sides of (\ref{section6.2-the proof of m_t-1 norm expectation}), we obtain
\begin{align}
\label{ling00}
\mathbb{E}&\bigl[\|\textbf{m}_{t-1}\|_F^2\bigr]
\le\theta_{t-1}\sum_{k=0}^{t-1}
{{\beta ^{t - 1 - k}}}\,
\mathbb{E}\!\left[\mathbb{E}_{\mathcal{L}_{t-1}}\bigl[\|\textbf{g}_k\|_F^2\bigr]\right]\nonumber\\
&=\theta_{t-1}\sum_{k=0}^{t-1}
{{\beta ^{t - 1 - k}}}
\mathbb{E}[\|{\textbf{g}}_k\|_F^2]\nonumber\\
&\le \theta_{t-1}\sum_{k=0}^{t-1}
{{\beta ^{t - 1 - k}}} \left(\frac{{8{d^2}\sigma _2^2 L_{{f_1}}^2}}{{\gamma _g^2}}\chi_k + \kappa _4\right)\nonumber\\
&=\theta_{t-1}\left(\frac{{8{d^2}\sigma _2^2 L_{{f_1}}^2}}{{\gamma _g^2}}\sum_{k=0}^{t-1}
{{\beta ^{t - 1 - k}}} \chi_k + \kappa _4\sum_{k=0}^{t-1}
{{\beta ^{t - 1 - k}}}\right)\nonumber\\
&\le \frac{1}{1-\beta}\left(\frac{{8{d^2}\sigma _2^2 L_{{f_1}}^2}}{{\gamma _g^2}}\sum_{k=0}^{t-1}
{{\beta ^{t - 1 - k}}} \chi_k + \kappa _4\sum_{k=0}^{t-1}
{{\beta ^{t - 1 - k}}}\right)\nonumber\\
&\le \frac{{8{d^2}\sigma _2^2 L_{{f_1}}^2}}{{(1 - \beta )\gamma _g^2}}\sum\limits_{k = 0}^{t - 1} {{\beta ^{t - 1 - k}}} \chi_k  + \frac{{{\kappa _4}}}{{(1 - \beta)^2 }},
\end{align}
where the second and last inequalities hold due to (\ref{0990}) and (\ref{111111}), respectively.

As is well known,
\begin{align}\label{section6.2-lemma 4-mt-1+gt-1 norm expectation-0}
    \mathbb{E} [\| {\eta (\beta {\mathbf{m}_{t - 1}} + {\mathbf{g}_{t - 1}})} \|_F^2] \le& 2{\eta ^2}\mathbb{E} [\| {\beta {\mathbf{m}_{t - 1}}}  \|_F^2] \nonumber\\
&+ 2{\eta ^2}\mathbb{E} [\| {{\mathbf{g}_{t - 1}}}  \|_F^2].
\end{align}
By combining (\ref{0990}) with (\ref{ling00}), (\ref{section6.2-lemma 4-mt-1+gt-1 norm expectation-0}) can be rearranged as (\ref{section3.2-lemma4-whole}).
\subsection{Auxiliary Results}\label{Auxiliary Results}
In this section, we provide the upper bound of $\sum\limits_{t=0}^{T - 1} {{\chi_t}}$. 

To facilitate the analysis, we define 
\begin{equation}\label{pop}
 \rho= \frac{{16{\eta ^2}{d^2}\sigma _2^2L_{{f_1}}^2}}{{\gamma _g^2}}.   
\end{equation}

From~(\ref{section6.2-lemma4-Xt express}),~(\ref{section6.2-lemma4-bounded of the expected norm squared}),~(\ref{ling00}) and (\ref{section6.2-lemma 4-mt-1+gt-1 norm expectation-0}), we have
\begin{align}
\label{section6.2-chit detail expression}
 \chi_t \le &\max \{ {\ell _1},{\ell _2}\} {\chi _{t - 1}} + {\ell _3}\left( {\frac{{\rho {\beta ^2}}}{{1 - \beta }}\sum\limits_{k = 0}^{t - 1} {{\beta ^{t - 1 - k}}} {\chi _k}} \right.\nonumber\\
&\left. { + \rho {\chi _{t - 1}} + 2{\eta ^2}{\kappa _4}(1 + \frac{{{\beta ^2}}}{{{{(1 - \beta )}^2}}})} \right)
\nonumber\\
= &(\max\{ {\ell _1},{\ell _2}\}+ {\rho \ell_3 }){\chi _{t - 1}}+ \frac{{\rho}{\ell _3}{{\beta ^2}}}{{1 - \beta }}\sum\limits_{k = 0}^{t - 1} {{\beta ^{t - 1 - k}}} {\chi _k} \nonumber\\
&+ 2{\ell _3}{\eta ^2}{\kappa _4}\left(1 + \frac{{{\beta ^2}}}{{{{(1 - \beta )}^2}}}\right).
\end{align}
Summing both sides of (\ref{section6.2-chit detail expression}) from $t=0$ to $T-1$, we have
\begin{align}
\label{section6.2chit sum}
\sum\limits_{t= 0}^{T - 1} {{\chi _t}} \le& (\max\{ {\ell _1},{\ell _2}\} + {\rho \ell_3 })\sum\limits_{t = 0}^{T - 1} {{\chi _{t - 1}}} \nonumber\\
&+ \frac{{\rho}{\ell _3}{{\beta ^2}}}{{1 - \beta}}\sum\limits_{t = 0}^{T - 1} {\sum\limits_{k = 0}^{t - 1} {{\beta ^{t - 1 - k}}} {\chi _k}} \nonumber\\
&+ 2T{\ell _3}{\eta ^2}{\kappa _4}\left(1 + \frac{{{\beta ^2}}}{{{{(1 - \beta )}^2}}}\right).
\end{align}
Then, we consider the term $\sum\limits_{t = 0}^{T - 1} {\sum\limits_{k = 0}^{t - 1} {{\beta ^{t - 1 - k}}} {\chi _k}}$ in (\ref{section6.2chit sum}).

More notably,
\begin{align}\label{6565}
\sum\limits_{t = 0}^{T - 1} {\sum\limits_{k = 0}^{t - 1} {{\beta ^{t - 1 - k}}} {\chi _k}}=\sum\limits_{t = 1}^{T - 1} {\sum\limits_{k = 0}^{t - 1} {{\beta ^{t - 1 - k}}} {\chi _k}}= \sum\limits_{j = 0}^{T - 2} {\sum\limits_{k = 0}^j {{\beta ^{j - k}}} {\chi _k}}.
\end{align}
Now, we prove
\begin{equation}\label{eq:double-sum-T}
\sum_{j=0}^{T-2}\sum_{k=0}^{j} \beta^{\,j-k}\chi_k=
\sum_{k=0}^{T-2}\sum_{j=k}^{T-2} \beta^{\,j-k}\chi_k.
\end{equation}
For $T=2$, we have
\[\sum_{j=0}^{0}\sum_{k=0}^{j} \beta^{j-k}\chi_k= \chi_0,\]
and
\[\sum_{k=0}^{0}\sum_{j=k}^{0} \beta^{j-k}\chi_k=\chi_0.\]
Hence, \eqref{eq:double-sum-T} holds for $T=2$.

Assume that \eqref{eq:double-sum-T} holds for $T=M\ge 2$, that is,
\begin{equation}\label{eq:IH-T}
\sum_{j=0}^{M-2}\sum_{k=0}^{j} \beta^{j-k}\chi_k=\sum_{k=0}^{M-2}\sum_{j=k}^{M-2} \beta^{j-k}\chi_k.
\end{equation}
We show that \eqref{eq:double-sum-T} also holds for $T=M+1$, that is,
\begin{align}\label{eq:goal-T+1}
  \sum_{j=0}^{M-1}\sum_{k=0}^{j} \beta^{j-k}\chi_k=\sum_{k=0}^{M-1}\sum_{j=k}^{M-1} \beta^{j-k}\chi_k.  
\end{align}
For the left-hand side of \eqref{eq:goal-T+1}, we have
\begin{align*}
\sum_{j=0}^{M-1}\sum_{k=0}^{j} &\beta^{j-k}\chi_k
=\sum_{j=0}^{M-2}\sum_{k=0}^{j} \beta^{j-k}\chi_k+ \sum_{k=0}^{M-1} \beta^{M-1-k}\chi_k\\
&=\sum_{k=0}^{M-2}\sum_{j=k}^{M-2} \beta^{j-k}\chi_k+ \sum_{k=0}^{M-1} \beta^{M-1-k}\chi_k\\
&= \sum\limits_{k = 0}^{M - 2} {\sum\limits_{j = k}^{M - 2} {{\beta ^{j - k}}} } {\chi _k} + \sum\limits_{k = 0}^{M - 2} {{\beta ^{M - 1 - k}}} {\chi _k} +{\chi _{M - 1}}\\
&=\sum_{k=0}^{M-2}\left(\sum_{j=k}^{M-2} \beta^{j-k}\chi_k+ \beta^{M-1-k}\chi_k\right)+\chi_{M-1} \\
&=\sum_{k=0}^{M-2}\sum_{j=k}^{M-1} \beta^{j-k}\chi_k+\chi_{M-1}\\
&=\sum_{k=0}^{M-1}\sum_{j=k}^{M-1} \beta^{j-k}\chi_k,
\end{align*}
which is equivalent to the right-hand side of~\eqref{eq:goal-T+1}. Therefore,~\eqref{eq:double-sum-T} holds for $T\ge 2$. 

From~(\ref{6565}) and (\ref{eq:double-sum-T}), we have
\begin{align}
\label{Mathematical induction}
\sum\limits_{t = 0}^{T - 1} {\sum\limits_{k = 0}^{t - 1} {{\beta ^{t - 1 - k}}} {\chi _k}}&= \sum\limits_{k = 0}^{T - 2} {\sum\limits_{j = k}^{T - 2} {{\beta ^{j - k}}} {\chi _k}}  \nonumber\\
&= \sum\limits_{k = 0}^{T - 2} {{\chi _k}\sum\limits_{s = 0}^{T - 2 - k} {{\beta ^s}} }\nonumber\\
&= \sum\limits_{k = 0}^{T - 2} {{\chi _k}} \frac{{1 - {\beta ^{T - 1 - k}}}}{{1 - \beta }}\nonumber\\
&\le \frac{1}{{1 - \beta }}\sum\limits_{k = 0}^{T - 2} {{\chi _k}}\nonumber\\
&\le \frac{1}{{1 - \beta }}\sum\limits_{k = 0}^{T - 1} {{\chi _k}}.
\end{align}
From~(\ref{section6.2chit sum}) and~(\ref{Mathematical induction}), we have
\begin{equation} 
\begin{aligned} \label{section6.2-chi_t sum simplified}
\sum\limits_{t = 0}^{T - 1} {{\chi _t}}  &\le (max\{ {\ell _1},{\ell _2}\}  + \rho {\ell _3})\sum\limits_{t = 0}^{T - 1} {{\chi _{t - 1}}}  + \frac{{\rho {\ell _3}{\beta ^2}}}{{{{(1 - \beta )}^2}}}\sum\limits_{t = 0}^{T - 1} {{\chi _t}}  \\
&~~~+ 2T{\ell _3}{\eta ^2}{\kappa _4}\left(1 + \frac{{{\beta ^2}}}{{{{(1 - \beta )}^2}}}\right).
\end{aligned}
\end{equation}
Based on $\sum\limits_{t = 0}^{T - 1} {{\chi _{t - 1}}}  = {\chi _{ - 1}} + \sum\limits_{t = 0}^{T - 2} {{\chi _t}} \le {\chi _{ - 1}} + \sum\limits_{t = 0}^{T - 1} {{\chi _t}}$, from (\ref{section6.2-chi_t sum simplified}), we have
\begin{align}
\label{1221}
\sum\limits_{t = 0}^{T - 1} {{\chi _t}}  \le& (\max\{ {\ell _1},{\ell _2}\}  + \rho {\ell _3})\left({\chi _{ - 1}} + \sum\limits_{t = 0}^{T - 1} {{\chi _t}}\right)\nonumber\\
&+ \frac{{\rho {\ell _3}{\beta ^2}}}{{{{(1 - \beta )}^2}}}\sum\limits_{t = 0}^{T - 1} {{\chi _t}}  + 2T{\ell _3}{\eta ^2}{\kappa _4}\left(1 + \frac{{{\beta ^2}}}{{{{(1 - \beta )}^2}}}\right)\nonumber\\
 =& (\max\{ {\ell _1},{\ell _2}\}  + \rho {\ell _3})\sum\limits_{t = 0}^{T - 1} {{\chi _t}}  + \frac{{\rho {\ell _3}{\beta ^2}}}{{{{(1 - \beta )}^2}}}\sum\limits_{t = 0}^{T - 1} {{\chi _t}} \nonumber\\
 &+ 2T{\ell _3}{\eta ^2}{\kappa _4}\left(1 + \frac{{{\beta ^2}}}{{{{(1 -\beta )}^2}}}\right)\nonumber\\
  =& \left(\max\{ {\ell _1},{\ell _2}\}  + \rho {\ell _3}(1+\frac{{{\beta ^2}}}{{{{(1 - \beta )}^2}}})\right)\sum\limits_{t = 0}^{T - 1} {{\chi _t}}\nonumber\\
 &+ 2T{\ell _3}{\eta ^2}{\kappa _4}\left(1 + \frac{{{\beta ^2}}}{{{{(1 -\beta )}^2}}}\right),  
\end{align}
where the first equality holds due to the settings of the Algorithm 1.
\subsection{Proof of Theorem 1}\label{appendix6.3-Proof of Theorem}
To simplify the analysis, we denote the following notations:
\begin{align}
{\tilde \varepsilon _0}&=\frac{{9{\lambda ^2}}}{{{\delta ^2}}} + \frac{4}{\omega } - 2,~~~~~~~~~~~
{\tilde \varepsilon _1}= \max \{ {c_1},...,{c_{12}}\},~\nonumber\\
{\tilde \varepsilon _2}&= 1500\phi \delta {\lambda ^2} + 25\phi \delta {\lambda ^2}{\omega ^2},~~~
{\tilde \varepsilon _3}=- \delta  + \frac{{45\phi {\lambda ^2}}}{{\omega \delta }},~\nonumber\\
\phi &= \frac{{2\rho }}{{{{(1 - \beta )}^2}}},~~
{\tilde \varepsilon _4}= {{\tilde \varepsilon }_0}\phi,~~{\tilde \varepsilon _5}= \max \{ {\tilde c_1},...,{\tilde c_{10}}\},\nonumber
\end{align}
 \begin{align}
{\tilde \varepsilon_6}&=\frac{{9{\lambda ^2}}}{2\delta } + \frac{{45\phi {\lambda ^2}}}{{\omega \delta }},~~~~
{\tilde \varepsilon _7}={{\tilde \varepsilon }_4}- \frac{\omega}{4},\nonumber\\
{{\tilde \varepsilon _8}} &=\frac{{27{\omega ^3}}}{{\delta (18{\delta ^2}{\omega ^3} + 9000{\lambda ^2} + 305{\delta ^2 \omega})}},\nonumber\\
{{\tilde \varepsilon }_9}&= \frac{{27{\omega ^2}}}{{\delta (1800{\lambda ^2} + 305{\delta ^2})}}, 
~~
{\phi _1} = n{d^2}\sigma _2^2L_{{f_1}}^2,\nonumber\\
{c_1} &=\frac{{640{\phi _1}{{\tilde \varepsilon }_0}}}{{\omega \gamma _g^2}},~~~
{c_2}=\frac{{2880{\phi _1}}}{{\omega \gamma _g^2}},~~~~~{c_3}=\frac{{54000{\phi _1}{\lambda ^2}}}{{{\delta ^2}{\omega ^2}\gamma _g^2}},\nonumber\\
{c_4} &=\frac{{48{\phi _1}}}{{{\delta ^3}\gamma _g^2}},~~~~~~~
{c_5}=\frac{{2880{\phi _1}{\lambda ^2}}}{{{\delta ^2}\omega \gamma _g^2}},~~
{c_6}=\frac{{128\tilde \varepsilon _0^2{\phi _1}}}{{\gamma _g^2}},\nonumber\\
{c_7} &=\frac{{5324000\phi _1}}{{\tilde \varepsilon _8^2{\delta ^6}\gamma _g^2}},~~
{c_8}=\frac{{60000{\phi _1}{\lambda ^4}}}{{\omega {\delta ^4}\gamma _g^2}},\nonumber\\
{c_9}&=\frac{{40000{\phi _1}}}{{{3\omega ^2}\gamma _g^2}},~~~~
{c_{10}}=\frac{{25n {d^2}L_{{f_2}}^2\sigma _1^2{\beta ^8}}}{{{{(1 - \beta )}^4}}},\nonumber\\
{c_{11}}&=\frac{{4{d^4}\sigma _1^2L_{{f_2}}^4 }}{n},~~~
{c_{12}}=\frac{{8{L_{{f_2}}}{d^2}\sigma _2^2L_{{f_1}}^2}}{{\gamma _g^2}},\nonumber
\end{align}
where ${\tilde c_1}$--${\tilde c_{10}}$ and $\rho$ can be found in Appendix~\ref{1234567} and (\ref{pop}), respectively.


More notably, we need to carefully choose the parameters $\alpha_1$--$\alpha_8$ in $\varepsilon_1$--$\kappa_3$ and $\ell_1$--$\ell_3$. For details,
\begin{align}\label{ppp}
&\alpha _1 = \frac{2}{{{\gamma _x}\delta }},~{\alpha _4}  = \frac{{{\gamma _x}\delta }}{2},~{\alpha _7} = \frac{4}{\omega },\nonumber\\
  {\alpha _\iota } &= \frac{\omega }{4},~\iota  = 2,3,5,6,8,~{m_2}< {\gamma _x} < {m_1},
\end{align}
where 
 \begin{align}
{{m_1}} =&~{\left(\frac{{3\phi {\omega ^3}}}{{2({\delta ^3}{\omega ^3} + 500\delta{\lambda ^2} + {{\tilde \varepsilon }_2})}}\right)^{\frac{1}{3}}},\nonumber\\
{{m_2}}=&~\frac{{{{\tilde \varepsilon }_4} + \sqrt {\tilde \varepsilon _4^2 + \frac{{27( - {{\tilde \varepsilon }_3})\phi }}{2\delta }} }}{{2( - {{\tilde \varepsilon }_3})}},\nonumber\\
{m_3} =&~{\left(\frac{{3{\omega ^3}\phi }}{{2(100\delta \omega {\lambda ^2} + {{\tilde \varepsilon }_2)}}}\right)^{\frac{1}{3}}},\nonumber\\
{m_{4,1}} =&~\frac{{ - {{\tilde \varepsilon }_7} - \sqrt {\tilde \varepsilon _7^2 - \frac{{27{{\tilde \varepsilon }_6}\phi }}{2\delta }} }}{{2{{\tilde \varepsilon }_6}}}, \nonumber\\
{m_{4,2}} =& \frac{{ - {{\tilde \varepsilon }_7} + \sqrt {\tilde \varepsilon _7^2 - \frac{{27{{\tilde \varepsilon }_6}\phi }}{2\delta }} }}{{2{{\tilde \varepsilon }_6}}},\nonumber
 \end{align}
where ${m_2} = \max \{ {m_2},{m_{4,1}}\}$ and ${m_1} = \min \{ {m_1},{m_3},{m_{4,2}}\}$.

We next prove that the following inequality holds under the parameter settings in (\ref{ppp}).
\begin{align}
\label{qqq}
&\max\{ {\ell _1},{\ell _2}\}  + \rho {\ell _3}\left(1 + \frac{{{\beta ^2}}}{{{{(1 - \beta )}^2}}}\right)\nonumber\\
&\le \max\{ {\ell _1},{\ell _2}\}  + \frac{{2{\rho \ell _3}}}{{{{(1 - \beta )}^2}}}\nonumber\\
&< 1, 
\end{align}
where the first inequality holds due to $1 + \frac{{{\beta ^2}}}{{{{(1 - \beta )}^2}}}\le1 + \frac{{1}}{{{{(1 - \beta )}^2}}} \le \frac{2}{{{{(1 - \beta )}^2}}}$. 

~~~~\textbf{(i)} In this part, we substitute~(\ref{ppp}) into the corresponding components of ${\ell_1}$--${\ell_3}$.

For ${\varepsilon _1}(1 + \alpha _1^{ - 1})$, based on ${\gamma _x}\delta < 1$, we have $(1 + \frac{{{\gamma _x}\delta }}{2})(1 - {\gamma _x}\delta ) = 1 - \frac{{{\gamma _x}\delta }}{2} - \frac{{\gamma _x^2{\delta ^2}}}{2} < 1 - \frac{{{\gamma _x}\delta }}{2}$. Then, 
\begin{align}\label{z1}
{\varepsilon _1}(1 + \alpha _1^{ - 1}) = {(1 + \frac{{{\gamma _x}\delta }}{2})^2}{(1 - {\gamma _x}\delta )^2}<{(1 - \frac{{{\gamma _x}\delta }}{2})^2}.
\end{align}
For ${\kappa _3}$, based on 
$1 + \frac{4}{\omega } \le \frac{5}{\omega }$, 
we have
\begin{align}\label{z2}
{\kappa _3} = (1 + \frac{4}{\omega })\gamma _x^2{\lambda ^2}(1 + \frac{4}{\omega })(1 + \frac{4}{\omega })\le \frac{{125}}{{{\omega ^3}}}\gamma _x^2{\lambda ^2}.
\end{align}
For $(1 + {\alpha _2})(1 - \omega )(1 + {\alpha _3})$, we have
\begin{align}\label{z4}
(1 + {\alpha _2})(1 - \omega )(1 + {\alpha _3}) &= (1 + \frac{\omega }{4})(1 - \omega )(1 + \frac{\omega }{4}) \nonumber\\
&\le 1 - \frac{\omega }{4}.
\end{align}
For ${\varepsilon _4}$, based on $1 + \frac{2}{{{\gamma _x}\delta }}<\frac{3}{{{\gamma _x}\delta }}$, $1 + \frac{{{\gamma _x}\delta }}{2}<\frac{3}{2}$ and (\ref{z4}),
we have
\begin{align} \label{z3}
{\varepsilon _4}& = (1 + \frac{2}{{{\gamma _x}\delta }})\gamma _x^2{\lambda ^2}(1 + \frac{{{\gamma _x}\delta }}{{2}})(1 + \frac{\omega }{4})(1 - \omega )(1 + \frac{\omega }{4})\nonumber\\
 &< \frac{{9{\gamma _x}{\lambda ^2}}}{2\delta }(1 - \frac{\omega }{4}).
\end{align}
For $\kappa _2$, based on ${(1 + \frac{\omega }{4})^3}(1 - \omega ) \le 1 - \frac{\omega }{4}$, we have
\begin{align}\label{z5}
{\kappa _2} &=(1 + \frac{4}{\omega })\gamma _x^2{\lambda ^2}(1 + \frac{4}{\omega })(1 + \frac{\omega }{4})(1 + \frac{\omega }{4})(1 - \omega )(1 + \frac{\omega }{4})\nonumber\\
& \le \frac{{25}}{{{\omega ^2}}}\gamma _x^2{\lambda ^2}(1 - \frac{\omega }{4}). 
\end{align}
For ${\varepsilon _3}$, based on $1 + \frac{\omega }{4} \le \frac{5}{4}$ and $(1 + \frac{\omega }{4})(1 - \omega )(1 + \frac{4}{\omega }) + (1 + \frac{4}{\omega }) \le \frac{{10}}{\omega }$, we have
\begin{align}
\label{z6}
{\varepsilon _3} = &(1 + \frac{2}{{{\gamma _x}\delta }})\gamma _x^2{\lambda ^2}(1 + \frac{{{\gamma _x}\delta }}{2})((1 + \frac{\omega }{4})(1 - \omega )(1 + \frac{4}{\omega })\nonumber \\
&+ (1 + \frac{4}{\omega }))
<
\frac{{45{\gamma _x}{\lambda ^2}}}{{\omega \delta }}.    
\end{align}
For $(1 + {\alpha _1})({\varepsilon _1} + {\varepsilon _2})$, based on $(1 + \frac{{{\gamma _x}\delta }}{2}){(1 - {\gamma _x}\delta )^2} < (1 - \frac{{{\gamma _x}\delta }}{2})(1 - {\gamma _x}\delta ) <1$, we have
\begin{align}\label{z7}
 (1 +{\alpha _1})({\varepsilon _1} + {\varepsilon _2}) =& (1 + \frac{2}{{{\gamma _x}\delta }})(1 + \frac{{{\gamma _x}\delta }}{2}){(1 - {\gamma _x}\delta )^2} \nonumber\\
 &+ {(1 + \frac{2}{{{\gamma _x}\delta }})^2}\gamma _x^2{\lambda ^2}\nonumber\\
 <&\frac{3}{{{\gamma _x}\delta }}(1 + \frac{{3{\gamma _x}{\lambda ^2}}}{\delta }).   
\end{align}
For ${\kappa _1}$, we have
\begin{align}\label{z8}
{\kappa _1} =& (1 + \frac{4}{\omega })\gamma _x^2{\lambda ^2}( (1 + \frac{4}{\omega })(1 + \frac{\omega }{4})((1 + \frac{\omega }{4})(1 - \omega )\nonumber\\
&(1 + \frac{4}{\omega }) + (1 + \frac{4}{\omega })) + (1 + \frac{\omega }{4}))\nonumber\\
\le & (\frac{{375}}{{{\omega ^3}}} + \frac{{25}}{{4\omega }})\gamma _x^2{\lambda ^2}.   
\end{align}    
For $(1 + {\alpha _2})(1 - \omega )(1 + {\alpha _3^{-1}})$, based on $(1 + \frac{\omega }{4})(1 - \omega ) \le 1 - \frac{{3\omega }}{4}$, we have
\begin{align}   
\label{z9} 
(1 + {\alpha _2})(1 - \omega )(1 + \alpha _3^{ - 1}) &= (1 + \frac{\omega }{4})(1 - \omega )(1 + \frac{4}{\omega })\nonumber\\ 
&\le (1 - \frac{{3\omega }}{4})(1 + \frac{4}{\omega }) \nonumber\\ 
&\le \frac{4}{\omega }-2.
\end{align}    
From (\ref{z1})--(\ref{z9}), we have
\begin{subequations}
\begin{align}
    {\ell _1}&= {\varepsilon _1}(1 + \alpha _1^{ - 1}) + {\kappa _3}<{(1 - \frac{{{\gamma _x}\delta }}{2})^2} + \frac{{125\gamma _x^2{\lambda ^2}}}{{{\omega ^3}}},\label{section6.2-ell1 with detail parameter}\\
{\ell _2}& = {\varepsilon _4} + (1 + {\alpha _2})(1 - \omega )(1 + {\alpha _3}) + {\kappa _2} \nonumber\\
&< (1 - \frac{\omega }{4})(\frac{{9{\gamma _x}{\lambda ^2}}}{2\delta } + 1 + \frac{{25}}{{{\omega ^2}}}\gamma _x^2{\lambda ^2})\nonumber\\
&<1 - \frac{\omega }{4}+\frac{{9{\gamma _x}{\lambda ^2}}}{2\delta } + \frac{{25}}{{{\omega ^2}}}\gamma _x^2{\lambda ^2},\label{section6.2-ell2 with detail parameter}\\
{\ell _3}& = {\varepsilon _3} + (1 + {\alpha _1})({\varepsilon _1} + {\varepsilon _2})+ (1 + {\alpha _2})(1 - \omega )(1 + {\alpha _3}) + {\kappa _1}\nonumber\\
 &< \frac{{45{\gamma _x}{\lambda ^2}}}{{\omega \delta }} + \frac{3}{{{\gamma _x}\delta }}(1 + \frac{{3{\gamma _x}{\lambda ^2}}}{\delta }) +\frac{ 4}{\omega}-2 \nonumber\\
 &~~~+ (\frac{{375}}{{{\omega ^3}}} + \frac{{25}}{{4\omega }})\gamma _x^2{\lambda ^2}.\label{section6.2-ell3 with detail parameter}
\end{align}
\end{subequations}
~~~~\textbf{(ii)} In this part, we first show that ${\sqrt {\tilde \varepsilon _7^2 - \frac{{27{{\tilde \varepsilon }_6}\phi }}{2\delta }} }>0$ in ${m_{4,1}}$ and ${m_{4,2}}$. We then bound ${m_1},{m_3}$ and ${m_{4,2}}$ to establish that ${m_1} = \min \{ {m_1},{m_3},{m_{4,2}}\}$.

First, we prove that ${\sqrt {\tilde \varepsilon _7^2 - \frac{{27{{\tilde \varepsilon }_6}\phi }}{2\delta }} }>0$. 

From $\phi  < \frac{\omega }{{20{{\tilde \varepsilon }_0}}}$, we have  
\begin{align}\label{6666}
- {{\tilde \varepsilon }_7} = \frac{\omega }{4} - {{\tilde \varepsilon }_0}\phi  > \frac{\omega }{4} - \frac{\omega }{20} = \frac{\omega }{5}.
\end{align}   
Then, we have $(- {{\tilde \varepsilon }_7})^2 > \frac{{{\omega ^2}}}{{25}}$.

From $\phi  < \frac{\omega }{{90}}$, we have
\begin{align}\label{2345}
{{\tilde \varepsilon }_6} = \frac{{9{\lambda ^2}}}{\delta }\left(\frac{1}{2} + \frac{{5\phi }}{\omega }\right) < \frac{{5{\lambda ^2}}}{\delta }.
\end{align}
From $\phi  < \frac{{2{\omega ^2}{\delta ^2}}}{{3375{\lambda ^2}}}$ and~(\ref{6666})--(\ref{2345}), we have 
\begin{align}\label{clj000}
\tilde \varepsilon _7^2 - \frac{{27{{\tilde \varepsilon }_6}\phi }}{2\delta } > \frac{{{\omega ^2}}}{{25}} - \frac{{135\phi {\lambda ^2}}}{{{2\delta ^2}}} > 0.    
\end{align}  
Second, we analyze the upper bound of $m_1$. 

It is easy to deduce that
\begin{align*}
{\delta ^3}{\omega ^3} + 500\delta {\lambda ^2} + {\tilde \varepsilon _2}>{\delta ^3}{\omega ^3}.        
\end{align*}
From $\phi <\frac{{{2\delta ^3}}}{3}$, we have
\begin{align}\label{099}
{m_1} < {\left(\frac{{3\phi }}{{2{\delta ^3}}}\right)^{\frac{1}{3}}} < 1.
\end{align}
Third, we analyze the lower bound of $m_1$. Based on $\omega \in (0,1]$, we have
\begin{align*}
{{\tilde \varepsilon }_2} \le 1525\phi \delta {\lambda ^2}.    
\end{align*} 
From $\phi<\frac{{{\delta ^2}\omega }}{{90{\lambda ^2}}}$, we have
\begin{align}\label{zzz}
{\delta ^3}{\omega ^3} + 500\delta {\lambda ^2} + {{\tilde \varepsilon }_2} 
&\le {\delta ^3}{\omega ^3} + 500\delta {\lambda ^2}+ 1525\phi \delta {\lambda ^2}\nonumber\\
&< \frac{{\delta}}{{18}}(18{\delta ^2}{\omega ^3} +9000{\lambda ^2} + 305{\delta ^2 \omega}).    
\end{align} 
 Then, 
 \begin{align}\label{edc}
{m_1}> {\left({{\tilde \varepsilon }_8}\phi \right)^{\frac{1}{3}}}.
 \end{align}
Fourth, we analyze the lower bound of $m_3$. Similar to (\ref{zzz}), we have
\begin{align*}
100\delta \omega {\lambda ^2} + {{\tilde \varepsilon }_2}\le&
100\delta \omega {\lambda ^2} + 1525\phi \delta {\lambda ^2} \nonumber\\
<&
\frac{{\delta \omega }}{{18}}(1800{\lambda ^2} + 305{\delta ^2}).    
\end{align*}
 Then, 
\begin{align}\label{edc1}
{m_3} > {\left({{\tilde \varepsilon }_9}\phi \right)^{\frac{1}{3}}}.
 \end{align}
Now, we compare $m_1$ and $m_3$. Based on the definition of ${{\tilde \varepsilon }_2}$, we have
\begin{align*}
\frac{1}{{{\delta ^3}{\omega ^3} + 500\delta {\lambda ^2} + {{\tilde \varepsilon }_2}}} <\frac{1}{100\delta {\lambda ^2} + {{\tilde \varepsilon }_2}} < \frac{1}{{100\delta \omega {\lambda ^2} + {{\tilde \varepsilon }_2}}}.  
\end{align*}
 Then, we have
\begin{align}\label{rfv}
{m_1} < {m_3}.
\end{align}

Finally, we analyze the lower bound of $m_{4,2}$. 

More notably,
\begin{align}\label{huhuhu}
{m_{4,2}} 
> \frac{{ - {{\tilde \varepsilon }_7}}}{{2{{\tilde \varepsilon }_6}}}>\frac{{\omega \delta }}{{50{\lambda ^2}}},
\end{align}
where the first inequality holds due to (\ref{clj000}); the second inequality holds due to~(\ref{6666}) and (\ref{2345}).

From $100\delta \omega {\lambda ^2} + {{\tilde \varepsilon }_2}>
100\delta \omega {\lambda ^2}$, we have
\begin{align*}
{m_3} < {\left( {\frac{{3{\omega ^2}\phi }}{{200\delta {\lambda ^2}}}} \right)^{\frac{1}{3}}}.    
\end{align*}
From $\phi <\frac{{\omega {\delta ^4}}}{{1875{\lambda ^4}}}$, we have
\begin{align}\label{7uj}
{m_3}<{m_{4,2}}.
\end{align}
From (\ref{rfv}) and (\ref{7uj}), we have $\min \{ {m_1},{m_3},{m_{4,2}}\}={m_1}$. 

~~~~\textbf{(iii)} Note that~(\ref{qqq}) holds if ${\ell _1} + \phi {\ell _3}<1$ and ${\ell _2}+\phi {\ell _3} <1$. 

First, we analyze ${\ell _1} + \phi {\ell _3} <1$. From~(\ref{section6.2-ell1 with detail parameter}) and (\ref{section6.2-ell3 with detail parameter}), we have
\begin{align}\label{5665}
{\ell _1} + \phi {\ell _3}< & {\tilde \varepsilon _3}{\gamma _x} 
+ \left(\frac{{{\delta ^2}}}{4} + \frac{{125{\lambda ^2}}}{{{\omega ^3}}} + \frac{{375\phi {\lambda ^2}}}{{{\omega ^3}}} + \frac{{25\phi {\lambda ^2}}}{{4\omega }}\right)\gamma _x^2 \nonumber\\
&+ 1 + \frac{{3\phi }}{{{\gamma _x}\delta }} + {\tilde \varepsilon _4}.  
\end{align}

From ${\gamma _x}< m_1$, we have
\begin{align*}
\left(\frac{{{\delta ^2}}}{4} + \frac{{125{\lambda ^2}}}{{{\omega ^3}}} + \frac{{375\phi {\lambda ^2}}}{{{\omega ^3}}} + \frac{{25\phi {\lambda ^2}}}{{4\omega }}\right)\gamma _x^2< \frac{{3\phi }}{{{8\gamma _x}\delta }}. 
\end{align*}
Hence, from (\ref{5665}), we have
\begin{align}\label{gfd}
{\ell _1} + \phi {\ell _3}< & 1 + {\tilde \varepsilon _3}{\gamma _x} 
+ \frac{{27\phi }}{{8{\gamma _x}\delta }} + {\tilde \varepsilon _4}.     
\end{align}
Then, from ${\gamma _x} > {m_2}>0$, we have
\begin{equation*}\label{section6.2-the first gamma_x}
{{\tilde \varepsilon }_3}\gamma _x^2 + {{\tilde \varepsilon }_4}{\gamma _x} + \frac{{27\phi }}{{8\delta }} < 0.
\end{equation*}
That is,
\begin{equation}
{{\tilde \varepsilon }_3}\gamma _x + {{\tilde \varepsilon }_4} + \frac{{27\phi }}{{{8\gamma _x}\delta }} < 0.
\end{equation}
Now, we analyze the upper bound of $m_2$.

From $\phi < \frac{{{\delta ^2}\omega }}{{90{\lambda ^2}}}$, we have
\begin{align}\label{wsx}
-{{\tilde \varepsilon}_3}=\delta  - \frac{{45\phi {\lambda ^2}}}{{\omega \delta }}> \frac{\delta}{2}. 
\end{align}
It can be easily deduced that
\begin{align}
\tilde \varepsilon _4^2 + \frac{{27( - {{\tilde \varepsilon }_3})\phi }}{2\delta } &= \tilde \varepsilon _0^2{\phi ^2} + \frac{{27\phi }}{2\delta }\left(\delta  - \frac{{45\phi {\lambda ^2}}}{{\omega \delta }}\right) \nonumber\\
&< \tilde \varepsilon _0^2{\phi ^2} + \frac{{27}}{2}\phi.    
\end{align}
From $\phi<\frac{{1}}{{4\tilde \varepsilon _0^2}}$, we have
\begin{align}\label{qaz}
\tilde \varepsilon _0^2{\phi ^2} + \frac{{27}}{2}\phi< \frac{55\phi}{4}.    
\end{align}
From~(\ref{wsx})--(\ref{qaz}), we have
\begin{align}\label{4343}
{m_2}<\frac{{\sqrt {\tilde \varepsilon _4^2 + \frac{{27( - {{\tilde \varepsilon }_3})\phi }}{2\delta }} }}{{( - {{\tilde \varepsilon }_3})}}< \frac{{\sqrt{55}}}{\delta }{\phi ^{\frac{1}{2}}},
\end{align}
where the first inequality holds due to $\tilde{\varepsilon}_4 < \sqrt{\tilde{\varepsilon}_4^2 + \frac{27(-\tilde{\varepsilon}_3)\phi}{2\delta}}.$
From (\ref{edc}),~(\ref{4343}) and $\phi < \frac{{\tilde \varepsilon _8^2{\delta ^6}}}{{166375}}$, we have 
\begin{align}\label{xxx}
{m_2} < \frac{{{\sqrt{55}}}}{\delta }{\phi ^{\frac{1}{2}}}< {({{\tilde \varepsilon }_8}\phi )^{\frac{1}{3}}} < {m_1}.    
\end{align}
From (\ref{gfd}), we have ${\ell _1} + \phi {\ell _3} < 1$ when $m_2<\gamma_x< m_1$.

Second, we analyze ${\ell _2} + \phi {\ell _3} <1$. 

Based on (\ref{section6.2-ell2 with detail parameter}) and (\ref{section6.2-ell3 with detail parameter}), we have
\begin{align}\label{1qaz}
{\ell _2} + \phi {\ell _3}<& {{\tilde \varepsilon _6}}{\gamma _x}+ \left(\frac{{25{\lambda ^2}}}{{{\omega ^2}}} + \frac{{375\phi {\lambda ^2}}}{{{\omega ^3}}} + \frac{{25\phi {\lambda ^2}}}{{4\omega }}\right)\gamma _x^2 + \frac{{3\phi }}{{{\gamma _x}\delta }}\nonumber \\
&+ 1 + {{\tilde \varepsilon _7}}.
\end{align}
From ${\gamma _x}< {m_3}$, we have
\begin{align*}
\left(\frac{{25{\lambda ^2}}}{{{\omega ^2}}} + \frac{{375\phi {\lambda ^2}}}{{{\omega ^3}}} + \frac{{25\phi {\lambda ^2}}}{{4\omega }}\right)\gamma _x^2<\frac{{3\phi }}{{{8\gamma _x}\delta }}.    
\end{align*}
Hence, from (\ref{1qaz}), we have 
\begin{align}
 {\ell _2} + \phi {\ell _3}<&1+ {{\tilde \varepsilon _6}}{\gamma _x}+ \frac{{27\phi }}{{{8\gamma _x}\delta }}\nonumber + {{\tilde \varepsilon _7}}.   
\end{align}
It is clear that 
\begin{align}\label{sososo}
{{\tilde \varepsilon }_6} = \frac{{9{\lambda ^2}}}{\delta }(\frac{1}{2} + \frac{{5\phi }}{\omega }) > \frac{{9{\lambda ^2}}}{2\delta }.    
\end{align}
When $0<{m_{4,1}} < {\gamma _x} < {m_{4,2}}$, we have
\begin{equation}\label{section6.2-the second gamma x}
 {{\tilde \varepsilon _6}}\gamma _x^2 + {{\tilde \varepsilon _7}}{\gamma _x} +\frac{{27\phi }}{{8\delta }} < 0.
\end{equation}
That is,
\begin{align*}
 {{\tilde \varepsilon _6}}\gamma _x+ {{\tilde \varepsilon _7}}+\frac{{27\phi }}{{8\gamma _x} {\delta }} < 0.    
\end{align*}
Hence, from (\ref{7uj}) and (\ref{section6.2-the second gamma x}), we have ${\ell _2} + \phi {\ell _3} < 1$ when $m_{{4,1}}<\gamma_x< m_3$.

Now, we bound $m_{4,1}$ and $m_2$ to analyze ${m_2} = \max \{ {m_2},{m_{4,1}}\}$.

First, we analyze the upper bound of $m_{4,1}$. 

From~(\ref{clj000}), (\ref{section6.2-the second gamma x}) and (\ref{sososo}), we have
\begin{align}
{m_{4,1}} \cdot {m_{4,2}} = \frac{{27\phi }}{{8{{\tilde \varepsilon }_6}\delta }} < \frac{{3\phi}}{{4{\lambda ^2}}}. 
\end{align}
Then, from (\ref{huhuhu}), we have
\begin{align}\label{2332}
{m_{4,1}} = \frac{{{m_{4,1}} \cdot {m_{4,2}}}}{{{m_{4,2}}}} < \frac{{75\phi }}{{2\omega \delta }}. 
\end{align}
Second, we analyze the lower bound of $m_2$.

From (\ref{wsx}), we have
\begin{align*}
 \tilde \varepsilon _4^2 + \frac{{27( - {{\tilde \varepsilon }_3})\phi }}{{2\delta }} > \tilde \varepsilon _0^2{\phi ^2} + \frac{{27}}{4}\phi. 
\end{align*}
Then, based on $-{{\tilde \varepsilon}_3}<{\delta}$, we have
\begin{align*}
{m_2} > \frac{{{{\tilde \varepsilon }_4} + \sqrt {\tilde \varepsilon _4^2 + \frac{{27}}{4}\phi } }}{{2( - {{\tilde \varepsilon }_3})}} > \frac{{\sqrt {\tilde \varepsilon _4^2 + \frac{{27}}{4}\phi } }}{{2( - {{\tilde \varepsilon }_3})}} > \frac{{\sqrt {\frac{{27}}{4}\phi } }}{{2( - {{\tilde \varepsilon }_3})}} > \frac{{3\sqrt {3\phi } }}{{2\sqrt {2}\delta }}. 
\end{align*}

From $\phi<\frac{{3{\omega ^2}}}{{1250}}$, we have
\begin{align}\label{2112}
{m_{4,1}} <{m_2}. 
\end{align}
It follows that (\ref{qqq}) holds when $m_2<\gamma_x< m_1$.
Hence, from (\ref{1221}), we have
\begin{align}
\label{section6.2-chi_t summary finally}
&\left( {1 - \max \{ {\ell _1},{\ell _2}\}  - \rho {\ell _3}(1 + \frac{{{\beta ^2}}}{{{{(1 - \beta )}^2}}})} \right)\sum\limits_{t = 0}^{T - 1} {{\chi _t}}  \nonumber \\
&\le 2T{\ell _3}{\eta ^2}{\kappa _4}\left(1 + \frac{{{\beta ^2}}}{{{{(1 - \beta )}^2}}}\right) \nonumber \\
&\le \frac{{{4T{\ell _3}{\eta ^2}{\kappa _4}}}}{{{{(1 - \beta) }^2}}}.
\end{align}
For the sake of analysis, denote $\rho_1=max\{ {\ell _1},{\ell _2}\}  + \rho {\ell _3}(1 + \frac{{{\beta ^2}}}{{{{(1 - \beta )}^2}}})$. From (\ref{section6.2-chi_t summary finally}), we have
\begin{align}\label{section6.2-lemma4-chit summary finally before set parameters}
\sum\limits_{t = 0}^{T - 1} {{\chi _t}}  
\le \frac{{4T{\ell _3}{\eta ^2}{\kappa _4}}}{{{\rm{(1}} - {\rho _1}){{(1 - \beta )}^2}}}.
\end{align}
By averaging over~(\ref{section6.2-lemma4-chit summary finally before set parameters}) over $t=0$ to $t=T-1$, we obtain 
\begin{align}\label{section6.2-lemma4-chi t iteartion results}
\frac{1}{T}\sum\limits_{t = 0}^{T - 1} {{\chi _t}}  \le \frac{{4{\ell _3}{\eta ^2}{\kappa _4}}}{{{\rm{(1}} - {\rho _1}){{(1 - \beta )}^2}}}.
\end{align}
Based on $\eta  = (1 - \beta )\sqrt {\frac{{n}}{T}}$
and $\frac{1}{T}\sum\limits_{t = 0}^{T - 1} {\sum\limits_{i = 1}^n {\mathbb{E}{{[\| {{x_{i,t}} - {{\bar x}_t}} \|}^2}} ]}
\le \frac{1}{T}\\\sum\limits_{t = 0}^{T - 1} {{\chi _t}}$. Then, from (\ref{light}) and (\ref{section6.2-lemma4-chi t iteartion results}), we have
\begin{align}\label{section3.2-theorem1-consensus error111}
 \frac{1}{T}\sum\limits_{t = 0}^{T - 1} &{\sum\limits_{i = 1}^n {\mathbb{E}{{[\| {{x_{i,t}} - {{\bar x}_t}} \|}^2]}} }\le\frac{{4{\ell _3}n{\kappa _4}}}{{{T(1-\rho_1)}}}\nonumber\\
\le& \frac{{16{n^2}{\ell _3}{d^2}\sigma _2^4L_{{f_1}}^2}}{{{{T(1-\rho_1)}}}}+\frac{{32{n^2}{\ell _3}{d^2}\sigma _2^2(\gamma _1 + {\vartheta _1})}}{{{{T(1-\rho_1)}}\gamma _g^2}}.
\end{align}
~~~~\textbf{(iv)} This part shows the upper bound of $\frac{1}{T}\sum\limits_{t = 0}^{T - 1}\mathbb{E}[\| {\nabla f({{\bar x}_t})\|{^2}}] $.

Based on the\textit{ Assumption \ref{section 2.1-ass:the objective function}} and (\ref{section6.2-lemma4-recurrence relation of the average of auxiliary sequence}), we have 
\begin{equation}\label{section6.3-theorem1-smooth of ecurrence relation of the average of vauxiliary sequence}
\begin{aligned}
     {f}({\bar x}_{t+1}^a) & \le f(\bar x_t^a) +  \langle {\nabla f(\bar x_t^a),\bar x_{t + 1}^a - \bar x_t^a}  \rangle  + \frac{{{L_{{f_2}}}}}{2}{ \| {\bar x_{t + 1}^a - \bar x_t^a} \|^2}\\
      &\le {f}({\bar x}_t^a) -\langle {\nabla {f}({\bar x}_t^a),\frac{\eta }{{1 - \beta }}{{\bar g}_t}} \rangle + \frac{{L_{{f_2}}}}{2}{(\frac{\eta }{{1 - \beta }})^2}{\| {{{\bar g}_t}} \|^2},
\end{aligned}
\end{equation}
where ${{\bar g}_t}={\frac{1}{n}\sum\limits_{i = 1}^n {{g_{i,t}}} }$. Based on the \textit{Lemma~\ref{section2.3-lemma1}}, we have
\begin{align}\label{section6.3-theorem1-gradient average expectation}
\mathbb{E}_{{{\cal L}_t}}[{{\bar g}_t}] &=\frac{1}{n}\sum\limits_{i = 1}^n {\mathbb{E}_{{{\cal L}_t}}[{g_{i,t}}]}\nonumber\\
&= d{\sigma _1}[\frac{1}{n}\sum\limits_{i = 1}^n {(\nabla {F_i}({x_{i,t}}) + {b_{i,t}})} ] \nonumber\\
    &= d{\sigma _1}[\nabla {{\bar F}}({x_{t}}) + {{\bar b}_{t}}].
\end{align}
where $\nabla {{\bar F}}({x_{t}}) = \frac{1}{n}\sum\limits_{i = 1}^n {\nabla {F_i}({x_{i,t}})}$. Then,
\begin{align}\label{section6.3-theorem1-biased norm}
   \|{b_{i,t}}\|& \le\frac{{{\gamma _g}}}{{2{\sigma _1}}}\mathbb{E}{_u}[ {{{\| {{u_{i,t}}}\|}_2}{{\| {{u_{i,t}}} \|}_2}{{\| {{\nabla ^2}{F_i}(\upsilon)} \|}_2}{{\| {{u_{i,t}}} \|}_2}} ]\nonumber\\
   &\le \frac{{{\gamma _g}}}{{2{\sigma _1}}}\sigma _2^3{L_{{f_2}}},
\end{align}
hence, $\|{\bar b_{t}}\|\le\frac{{{\gamma _g}}}{{2{\sigma _1}}}\sigma _2^3{L_{{f_2}}}$.
Thus, 
\begin{align*}
    \| {\nabla {f}({\bar x}_t^a)} - {\nabla {{\bar F}}({x_{t}})} \| =& \left\| {\frac{1}{n}\sum\limits_{i = 1}^n {(\nabla {F_i}({\bar x}_t^a) - \nabla {F_i}({x_{i,t}}))} } \right\| \nonumber
    \end{align*}
    \begin{align}
 \label{section6.3-theorem1-two gradient cha}
\le &   \frac{{L_{{f_2}}}}{n}\sum\limits_{i = 1}^n {\| {{\bar x}_t^a - {x_{i,t}}}\|}\nonumber\\
\le &  \frac{L_{{f_2}}}{{\sqrt n }}{( {\sum\limits_{i = 1}^n {{{\| {{\bar x}_t^a - {x_{i,t}}} \|}^2}} } )^{\frac{1}{2}}}, 
\end{align}
where the first inequality holds due to \textit{Assumption \ref{section 2.1-ass:the objective function}}; the second inequality holds due to the Cauchy--Schwarz inequality.

Taking the conditional expectation given ${{\cal L}_t}$ on both sides of (\ref{section6.3-theorem1-smooth of ecurrence relation of the average of vauxiliary sequence}), we have
\begin{align*}
\mathbb{E}&_{{{\cal L}_t}}[ {{f}}{({\bar x}_{t+1}^a)}]\nonumber\\
 \le& {f}({\bar x}_t^a)- \frac{\eta }{{1 - \beta }}\mathbb{E}_{{{\cal L}_t}}[ {\langle {\nabla {f}({\bar x}_t^a),{{\bar g}_t}} \rangle } ]\nonumber\\
 &+ \frac{{L_{{f_2}}}}{2}{( {\frac{\eta }{{1 - \beta }}} )^2}\mathbb{E}_{{{\cal L}_t}} {{{[\| {{{\bar g}_t}} \|}^2]}}\nonumber\\
\mathop  \le \limits^{(a)}   & {f}({\bar x}_t^a)- \frac{{\eta d{\sigma _1}}}{{1 - \beta }} {\langle {\nabla {f}({\bar x}_t^a),\nabla {{\bar F}}({x_{t}}) + {{\bar b}_t}} \rangle } \nonumber\\
&+ \frac{{L_{{f_2}}}}{2}{( {\frac{\eta }{{1 - \beta }}} )^2}\left(\frac{{8{d^2}\sigma _2^2 L_{{f_1}}^2}}{{n\gamma _g^2}}\sum\limits_{i = 1}^n \|{{x_{i,t}} - {{\bar x}_t}\|{^2}}  + \frac{{{\kappa _4}}}{n}\right)\nonumber\\
=& {f}({\bar x}_t^a) - \frac{{\eta d{\sigma _1}}}{{1 - \beta }}\langle \nabla f(\bar x_t^a),\nabla {{\bar F}}({x_{t}}) + {{\bar b}_t}{+ \nabla f(\bar x_t^a)}- {\nabla f(\bar x_t^a)}\rangle \nonumber\\
&+  \frac{{L_{{f_2}}}}{2}{( {\frac{\eta }{{1 - \beta }}} )^2}\left(\frac{{8{d^2}\sigma _2^2 L_{{f_1}}^2}}{{n\gamma _g^2}}\sum\limits_{i = 1}^n \|{{x_{i,t}} - {{\bar x}_t}\|{^2}}+ \frac{{{\kappa _4}}}{n}\right)\nonumber\\
 \le & {f}({\bar x}_t^a) - \frac{{\eta d{\sigma _1}}}{{1 - \beta }}{\| {\nabla f(\bar x_t^a)} \|^2} \nonumber\\
 &- \frac{{\eta d{\sigma _1}}}{{1 - \beta }}\langle {\nabla f(\bar x_t^a),{{\bar b}_t}} \rangle \nonumber\\
 &+ \frac{{\eta d{\sigma _1}}}{{1 - \beta }}\langle {\nabla f(\bar x_t^a),\nabla f(\bar x_t^a) - \nabla {{\bar F}}({x_{t}})} \rangle \nonumber\\
 &+  \frac{{L_{{f_2}} }}{2}{( {\frac{\eta }{{1 - \beta }}} )^2}\left(\frac{{8{d^2}\sigma _2^2 L_{{f_1}}^2}}{{n\gamma _g^2}}\sum\limits_{i = 1}^n \|{{x_{i,t}} - {{\bar x}_t}\|{^2}} + \frac{{{\kappa _4}}}{n}\right) \nonumber\\
 \mathop  \le \limits^{(b)} &{f}({\bar x}_t^a)- \frac{{\eta {d}{\sigma _1}}}{{1 - \beta }}{\| {\nabla {f}({\bar x}_t^a)} \|^2} \nonumber\\
 &+ \frac{{\eta {d}{\sigma _1}}}{{1 - \beta }}\| {\nabla {f}({\bar x}_t^a)} \|\| {{{\bar b}_t}} \| \nonumber\\
&+ \frac{{\eta {d}{\sigma _1}}}{{1 - \beta }}\| {\nabla {f}({\bar x}_t^a)} \|\| {\nabla {f}({\bar x}_t^a) - \nabla {{\bar F}}({x_{t}})}\| \nonumber\\
&+  \frac{{L_{{f_2}} }}{2}{( {\frac{\eta }{{1 - \beta }}} )^2}\left(\frac{{8{d^2}\sigma _2^2 L_{{f_1}}^2}}{{n\gamma _g^2}}\sum\limits_{i = 1}^n \|{{x_{i,t}} - {{\bar x}_t}\|{^2}}+ \frac{{{\kappa _4}}}{n}\right) \nonumber\\
\mathop  \le \limits^{(c)} & {f}({\bar x}_t^a) - \frac{{\eta {d}{\sigma _1}}}{{1 - \beta }}{\| {\nabla f(\bar x_t^a)} \|^2}\nonumber\\
&+ \frac{{\eta {d}{\sigma _1}}}{{1 - \beta }}(\frac{1}{4}{\| {\nabla f(\bar x_t^a)}\|^2} + {\| {{{\bar b}_t}} \|^2}) \nonumber\\
&+ \frac{{\eta {d}{\sigma _1}}}{{1 - \beta }}(\frac{1}{4}{\| {\nabla f(\bar x_t^a)} \|^2} + {\| {\nabla f(\bar x_t^a) - \nabla {{\bar F}}{{({x_{t}})}}} \|^2})\nonumber\\
&+ \frac{{L_{{f_2}} }}{2}{( {\frac{\eta }{{1 - \beta }}} )^2}\left(\frac{{8{d^2}\sigma _2^2 L_{{f_1}}^2}}{{n\gamma _g^2}}\sum\limits_{i = 1}^n \|{{x_{i,t}} - {{\bar x}_t}\|{^2}}+ \frac{{{\kappa _4}}}{n}\right)\nonumber\\
 =& {f}({\bar x}_t^a) - \frac{{\eta {d}{\sigma _1}}}{{2(1 - \beta )}} {\| {\nabla f(\bar x_t^a)}\|^2}+ \frac{{\eta {d}{\sigma _1}}}{{1 - \beta }}{\| {{{\bar b}_t}} \|^2}\nonumber\\
  &+ \frac{{\eta {d}{\sigma _1}}}{{1 - \beta }}{\| {\nabla f(\bar x_t^a) - \nabla {{\bar F}}{{({x_{t}})}}} \|^2} \nonumber
   \end{align*}
 \begin{align}
  \label{section6.3-theorem1-Taking the conditional expectation given Ht on both sides of f barx}
  &+ \frac{{L_{{f_2}} }}{2}{( {\frac{\eta }{{1 - \beta }}} )^2}\left(\frac{{8{d^2}\sigma _2^2 L_{{f_1}}^2}}{{n\gamma _g^2}}\sum\limits_{i = 1}^n \|{{x_{i,t}} - {{\bar x}_t}\|{^2}}+ \frac{{{\kappa _4}}}{n}\right)\nonumber\\
 \mathop  \le \limits^{(d)}  &{f}({\bar x}_t^a) - \frac{{\eta {d}{\sigma _1}}}{{2(1 - \beta )}}{\| {\nabla {f}({\bar x}_t^a)} \|^2} + \frac{{\eta d\gamma _g^2\sigma _2^6L_{{f_2}}^2}}{{4{\sigma _1}(1 - \beta )}}\nonumber\\
&{\rm{ + }}\frac{{\eta d{\sigma _1}}}{{1 - \beta }}\frac{{{L_{f_2}^2}}}{n}\sum\limits_{i = 1}^n {{{\| {{\bar x}_t^a - {x_{i,t}}} \|}^2}}  \nonumber\\
&+ \frac{{L_{{f_2}} }}{2}{( {\frac{\eta }{{1 - \beta }}} )^2}\left(\frac{{8{d^2}\sigma _2^2 L_{{f_1}}^2}}{{n\gamma _g^2}}\sum\limits_{i = 1}^n \|{{x_{i,t}} - {{\bar x}_t}\|{^2}}+ \frac{{{\kappa _4}}}{n}\right)\nonumber\\
 \mathop  \le \limits^{(e)}   &{f}({\bar x}_t^a) - \frac{{\eta d{\sigma _1}}}{{4(1 - \beta )}}{\| {\nabla {f}({{\bar x}_t})} \|^2} + \frac{{\eta d\gamma _g^2\sigma _2^6L_{{f_2}}^2}}{{4{\sigma _1}(1 - \beta )}}\nonumber\\
  &+ \frac{{\eta d {{{{L_{f_2}^2}}}}{\sigma _1}}}{{2(1 - \beta )}}{\| {{{\bar x}_t} - {\bar x}_t^a} \|^2}+\frac{{\eta d{\sigma _1}{{{L_{f_2}^2}}}}}{{{n}(1 - \beta) }} \sum\limits_{i = 1}^n {{{\| {{\bar x}_t^a - {x_{i,t}}} \|}^2}}\nonumber\\
  &+\frac{{4{L_{{f_2}}}{\eta ^2}{d^2}\sigma _2^2L_{{f_1}}^2}}{{n\gamma _g^2{{(1 - \beta )}^2}}}\sum\limits_{i = 1}^n \|{{x_{i,t}} - {{\bar x}_t}\|{^2}}+ \frac{{{\kappa _4}{L_{{f_2}}}{\eta ^2}}}{{2n{{(1-\beta)}^2}}}, 
\end{align}
where $(a)$ holds due to (\ref{section6.2-lemma4-gt average norm square and take expectation}) and (\ref{section6.3-theorem1-gradient average expectation}); $(b)$ and $(c)$ holds due to the Cauchy-Schwarz inequality and Young inequality, respectively. 
$(d)$ holds due to (\ref{section6.3-theorem1-biased norm}); $(e)$ holds due to ${\| {\nabla {f}({{\bar x}_t})}\|^2} \le 2{\| {\nabla {f}({{\bar x}_t}) - \nabla {f}({\bar x}_t^a)} \|^2} + 2{\| {\nabla {f}({\bar x}_t^a)}\|^2} \le 2{L_{f_2}^2}{\| {{{\bar x}_t} - {\bar x}_t^a}\|^2} + 2{\| {\nabla {f}({\bar x}_t^a)} \|^2}$. Rearranging this gives ${\| {\nabla {f}({\bar x}_t^a)}\|^2} \ge \frac{1}{2}{\| {\nabla {f}({{\bar x}_t})} \|^2} - {L_{f_2}^2}{\| {{{\bar x}_t} - {\bar x}_t^a} \|^2}$.

From~(\ref{section6.3-theorem1-Taking the conditional expectation given Ht on both sides
of f barx}), we have
\begin{align}
\label{T-14}
&\frac{{\eta d{\sigma _1}}}{{4(1 - \beta )}}{\| {\nabla f({{\bar x}_t})} \|^2}\nonumber\\
&\le f(\bar x_t^a) - \mathbb{E}_{{{\cal L}_t}}[f(\bar x_{t + 1}^a)] + \frac{{\eta d\gamma _g^2\sigma _2^6L_{{f_2}}^2}}{{4{\sigma _1}(1 - \beta )}}\nonumber\\
&~~~+ \frac{{\eta dL_{{f_2}}^2{\sigma _1}}}{{2(1 - \beta )}}{ \| {{{\bar x}_t} - \bar x_t^a} \|^2}\nonumber\\
 &~~~+ \frac{{\eta d{\sigma _1}L_{{f_2}}^2}}{{n(1 - \beta )}}\left(2\sum\limits_{i = 1}^n {{{ \| {{x_{i,t}} - {{\bar x}_t}}  \|}^2}}  + 2\sum\limits_{i = 1}^n {{{ \| {{{\bar x}_t} - \bar x_t^a}  \|}^2}} \right)\nonumber\\
&~~~+ \frac{{4{L_{{f_2}}}{\eta ^2}{d^2}\sigma _2^2L_{{f_1}}^2}}{{n\gamma _g^2{{(1 - \beta )}^2}}}{\sum\limits_{i = 1}^n \|{{x_{i,t}} - {{\bar x}_t}\|{^2}}} + \frac{{{\kappa _4}{L_{{f_2}}}{\eta ^2}}}{{2n{{(1 - \beta )}^2}}}\nonumber\\
 &\le f(\bar x_t^a) - \mathbb{E}_{{{\cal L}_t}}[f(\bar x_{t + 1}^a)] + \frac{{\eta d\gamma _g^2\sigma _2^6L_{{f_2}}^2}}{{4{\sigma _1}(1 - \beta )}} \nonumber\\
 &~~~+ \frac{{5\eta dL_{{f_2}}^2{\sigma _1}}}{{2(1 - \beta )}}{ \| {{{\bar x}_t} - \bar x_t^a}  \|^2}  \nonumber\\
 &~~~+\left(\frac{{4{L_{{f_2}}}{\eta ^2}{d^2}\sigma _2^2L_{{f_1}}^2}}{{n\gamma _g^2{{(1 - \beta )}^2}}}+\frac{{2\eta d{\sigma _1}L_{{f_2}}^2}}{{n(1 - \beta )}}\right){\sum\limits_{i = 1}^n \|{{x_{i,t}} - {{\bar x}_t}\|{^2}}}\nonumber\\
 &~~~+ \frac{{{\kappa _4}{L_{{f_2}}}{\eta ^2}}}{{2n{{(1 - \beta )}^2}}} .
\end{align}
From (\ref{section6.2-lemma4-the average of vauxiliary sequence}) and $\theta_{t-1}=\sum_{k=0}^{t-1}\beta^{t-1-k}=\frac{1-\beta^{t}}{1-\beta}$, we have 
\begin{align*}
\|\bar{x}_t - \bar{x}_t^{a}\|^{2}&=\left(\frac{\eta \beta^{2}}{1-\beta} \right)^{2}\left\| \frac{1}{n} \sum_{i=1}^{n}m_{i,t-1} \right\|^{2}\nonumber\\
&= \left(\frac{\eta \beta^{2}}{1-\beta} \right)^{2}\theta _{t - 1}^2\left\|\sum_{k=0}^{t-1} \frac{\beta^{t-1-k}}{\theta_{t-1}} \frac{1}{n} \sum_{i=1}^{n} g_{i,k}\right\|^{2}\nonumber
  \end{align*}
   \begin{align}\label{T-16}
&\le \left(\frac{\eta \beta^{2}}{1-\beta} \right)^{2}\theta _{t - 1}^2\sum_{k=0}^{t-1} \frac{\beta^{t-1-k}}{\theta_{t-1}}\left\| \frac{1}{n} \sum_{i=1}^{n} g_{i,k}\right\|^{2}\nonumber\\
&= \left(\frac{\eta \beta^{2}}{1-\beta} \right)^{2} \theta_{t-1}\sum_{k=0}^{t-1}{\beta^{t-1-k}}\left\| \frac{1}{n} \sum_{i=1}^{n} g_{i,k}\right\|^{2}\nonumber\\
 &\le  \frac{\beta^{4}\eta^{2}}{(1-\beta)^{3}}\sum_{k=0}^{t-1}{\beta^{t-1-k}}\left\| \frac{1}{n} \sum_{i=1}^{n} g_{i,k}\right\|^{2}.
\end{align}
From (\ref{T-14})--(\ref{T-16}), we have
\begin{align}\label{section6.2-theorem1-no expectation for gradient norm square} 
&\frac{{\eta d{\sigma _1}}}{{4(1 - \beta )}}{ \| {\nabla f({{\bar x}_t})} \|^2}\nonumber \\
&\le f(\bar x_t^a) - \mathbb{E}_{{{\cal L}_t}}[f(\bar x_{t + 1}^a)] + \frac{{\eta d\gamma _g^2\sigma _2^6L_{{f_2}}^2}}{{4{\sigma _1}(1 - \beta )}}+ \frac{{{\kappa _4}{L_{{f_2}}}{\eta ^2}}}{{2n{{(1 - \beta )}^2}}} \nonumber \\
&~~~+\frac{{5{\eta ^3}dL_{{f_2}}^2{\sigma _1}{\beta ^4}}}{{2{{(1-\beta )}^4}}}\sum\limits_{k = 0}^{t-1} {{\beta ^{t - 1 - k}}} {\left\| {\frac{1}{n}\sum\limits_{i = 1}^n {{g_{i,k}}} }  \right\|^2}\nonumber\\
&~~~+\left(\frac{{4{L_{{f_2}}}{\eta ^2}{d^2}\sigma _2^2L_{{f_1}}^2}}{{n\gamma _g^2{{(1 - \beta )}^2}}}+ \frac{{2\eta d{\sigma _1}L_{{f_2}}^2}}{{n(1 - \beta )}}\right){\sum\limits_{i = 1}^n \|{{x_{i,t}} - {{\bar x}_t}\|{^2}}}. 
\end{align}
Then, from~(\ref{section6.2-theorem1-no expectation for gradient norm square}), taking expectation with respect to the entire past and average over $t=0$ to $T-1$ gives 
\begin{align*}
&\frac{{\eta d{\sigma _1}}}{{4T(1 - \beta )}}\sum\limits_{t = 0}^{T - 1}  { \mathbb{E}[\| {\nabla f({{\bar x}_t})}  \|^2]}\nonumber\\
&\le \frac{f({{\bar x}_0}) - \mathbb{E}[f({{\bar x}_T^a})]}{T} + \frac{{\eta d\gamma _g^2\sigma _2^6L_{{f_2}}^2}}{{4{\sigma _1}(1 - \beta )}} + \frac{{{\kappa _4}{L_{{f_2}}}{\eta ^2}}}{{2n{{(1 - \beta )}^2}}}\nonumber\\
&~~~+ \frac{{5{\eta ^3}dL_{{f_2}}^2{\sigma _1}{\beta ^4}}}{{2{{(1 - \beta )}^4}}}\frac{1}{T}\sum\limits_{t = 0}^{T - 1} {\sum\limits_{k = 0}^{t - 1} {{\beta ^{t - 1 - k}}\left(\frac{{8{d^2}\sigma _2^2L_{{f_1}}^2}}{{n\gamma _g^2}}{\chi _k} + \frac{{{\kappa _4}}}{n}\right)} }\nonumber\\
&~~~+ \left(\frac{{4{L_{{f_2}}}{\eta ^2}{d^2}\sigma _2^2L_{{f_1}}^2}}{{n\gamma _g^2{{(1 - \beta )}^2}}} + \frac{{2\eta d{\sigma _1}L_{{f_2}}^2}}{{n(1 - \beta )}}\right)\frac{1}{T}\sum\limits_{t = 0}^{T - 1} {\chi _t} \nonumber\\
&=\frac{f({{\bar x}_0}) - \mathbb{E}[f({{\bar x}_T^a})]}{T} + \frac{{\eta d\gamma _g^2\sigma _2^6L_{{f_2}}^2}}{{4{\sigma _1}(1 - \beta )}} + \frac{{{\kappa _4}{L_{{f_2}}}{\eta ^2}}}{{2n{{(1 - \beta )}^2}}}\nonumber\\
&~~~+ \frac{{5{\eta ^3}dL_{{f_2}}^2{\sigma _1}{\beta ^4}}}{{2{{(1 - \beta )}^4}}}\left( {\frac{{8{d^2}\sigma _2^2L_{{f_1}}^2}}{{nT\gamma _g^2}}\sum\limits_{t = 0}^{T - 1} {\sum\limits_{k = 0}^{t - 1} {{\beta ^{t - 1 - k}}}{\chi _k} } } \right.\nonumber\\
&~~~\left. { + \frac{{{\kappa _4}}}{{nT}}\sum\limits_{t = 0}^{T - 1} {\sum\limits_{k = 0}^{t - 1} {{\beta ^{t - 1 - k}}} } } \right)\nonumber\\
&~~~+ \left(\frac{{4{L_{{f_2}}}{\eta ^2}{d^2}\sigma _2^2L_{{f_1}}^2}}{{n\gamma _g^2{{(1 - \beta )}^2}}} + \frac{{2\eta d{\sigma _1}L_{{f_2}}^2}}{{n(1 - \beta )}}\right)\frac{1}{T}\sum\limits_{t = 0}^{T - 1} {\chi _t} \nonumber\\
&\le\frac{f({{\bar x}_0}) - \mathbb{E}[f({{\bar x}_T^a})]}{T} + \frac{{\eta d\gamma _g^2\sigma _2^6L_{{f_2}}^2}}{{4{\sigma _1}(1 - \beta )}} + \frac{{{\kappa _4}{L_{{f_2}}}{\eta ^2}}}{{2n{{(1 - \beta )}^2}}}\nonumber\\
&~~~+ \frac{{5{\eta ^3}dL_{{f_2}}^2{\sigma _1}{\beta ^4}}}{{2{{(1 - \beta )}^4}}}\left( {\frac{{8{d^2}\sigma _2^2L_{{f_1}}^2}}{{nT\gamma _g^2({1 - \beta })}}\sum\limits_{k = 0}^{T - 1} {{\chi _k}} } \right.\left. { + \frac{{{\kappa _4}}}{{n({1 - \beta })}}} \right)\nonumber\\
&~~~+ \left(\frac{{4{L_{{f_2}}}{\eta ^2}{d^2}\sigma _2^2L_{{f_1}}^2}}{{n\gamma _g^2{{(1 - \beta )}^2}}} + \frac{{2\eta d{\sigma _1}L_{{f_2}}^2}}{{n(1 - \beta )}}\right)\frac{1}{T}\sum\limits_{t = 0}^{T - 1} {\chi _t} \nonumber\\
& \le \frac{f({{\bar x}_0}) - \mathbb{E}[f({{\bar x}_T^a})]}{T} + \frac{{\eta d\gamma _g^2\sigma _2^6L_{{f_2}}^2}}{{4{\sigma _1}(1 - \beta )}} + \frac{{{\kappa _4}{L_{{f_2}}}{\eta ^2}}}{{2n{{(1 - \beta )}^2}}} \nonumber
 \end{align*}
 \begin{align}
 \label{T-15}
 &~~~+ \left( {\frac{{4{L_{{f_2}}}{\eta ^2}{d^2}\sigma _2^2L_{{f_1}}^2}}{{n\gamma _g^2{{(1 - \beta )}^2}}} + \frac{{2\eta d{\sigma _1}L_{{f_2}}^2}}{{n(1 - \beta )}}} \right.\nonumber\\
&~~~\left. { + \frac{{20{\eta ^3}{d^3}L_{{f_2}}^2{\sigma _1}{\beta ^4}\sigma _2^2L_{{f_1}}^2}}{{{{(1 - \beta )}^5}n\gamma _g^2}}} \right) \frac{1}{T}\sum\limits_{t = 0}^{T - 1} {{\chi _t}}  + \frac{{5{\eta ^3}dL_{{f_2}}^2{\sigma _1}{\beta ^4}{\kappa _4}}}{{2n{{(1 - \beta )}^5}}},
\end{align}
where the first and second inequalities hold due to~(\ref{section6.2-lemma4-gt average norm square and take expectation}) and (\ref{Mathematical induction}), respectively. When $\eta  \le min\{\frac{{{{(1 - \beta )}^3}}}{{5d{L_{{f_2}}}{\sigma _1}{\beta ^4}}},\frac{{n(1 - \beta )}}{{2d{\sigma _1}L_{{f_2}}^2}}\}$, we have
\begin{align}
\label{section6.3-the gradient norm no chu coefficient}
&\frac{{\eta d{\sigma _1}}}{{4(1-\beta )}}\frac{1}{T}\sum\limits_{t = 0}^{T - 1} \mathbb{E}{{{ [\| {\nabla f({{\bar x}_t})}  \|}^2]}}\nonumber \\
&\le \frac{f({{\bar x}_0}) - \mathbb{E}[f({{\bar x}_T^a})]}{T}+ \frac{{\eta d\gamma _g^2\sigma _2^6L_{{f_2}}^2}}{{4{\sigma _1}(1 - \beta )}} + \frac{{{\kappa _4}{L_{{f_2}}}{\eta ^2}}}{{n{{(1 - \beta )}^2}}}\nonumber\\
&~~~+ \left(\frac{{8{L_{{f_2}}}{\eta ^2}{d^2}\sigma _2^2L_{{f_1}}^2}}{{n\gamma _g^2{{(1 - \beta )}^2}}} + 1\right)\frac{1}{T}\sum\limits_{t = 0}^{T - 1} {{\chi _t}}.
\end{align}
We further require the step size $\eta  \le \sqrt {\frac{{n\gamma _g^2{{(1 - \beta )}^2}}}{{8{L_{{f_2}}}{d^2}\sigma _2^2L_{{f_1}}^2}}}$ to establish
\begin{align}
\label{lololo}
&\frac{{\eta d{\sigma _1}}}{{4(1-\beta )}}\frac{1}{T}\sum\limits_{t = 0}^{T - 1} \mathbb{E}{{{ [\| {\nabla f({{\bar x}_t})}  \|}^2]}}\nonumber \\
&\le \frac{f({{\bar x}_0}) - \mathbb{E}[f({{\bar x}_T^a})]}{T} + \frac{{\eta d\gamma _g^2\sigma _2^6L_{{f_2}}^2}}{{4{\sigma _1}(1 - \beta )}} + \frac{{{\kappa _4}{L_{{f_2}}}{\eta ^2}}}{{n{{(1 - \beta )}^2}}}
+ \frac{2}{T}\sum\limits_{t = 0}^{T - 1} {{\chi _t}}\nonumber\\
&\le\frac{f({{\bar x}_0}) - \mathbb{E}[f({{\bar x}_T^a})]}{T} + \frac{{\eta d\gamma _g^2\sigma _2^6L_{{f_2}}^2}}{{4{\sigma _1}(1 - \beta )}}+ \frac{{{\kappa _4}{L_{{f_2}}}{\eta ^2}}}{{n{{(1 - \beta )}^2}}} \nonumber\\
&~~~~+ \frac{{8{\ell _3}{\eta ^2}{\kappa _4}}}{{{\rm{(1}} - {\rho _1}){{(1 - \beta )}^2}}}, 
\end{align}
where the second inequality holds due to~(\ref{section6.2-lemma4-chi t iteartion results}). 

Since $\mathbb{E}[f({{\bar x}_T^a})]\ge f^*$,
from~(\ref{light}) and (\ref{lololo}), we have
\begin{align}
\label{888}
\frac{1}{T}&\sum\limits_{t = 0}^{T - 1} \mathbb{E}{{{ [\| {\nabla f({{\bar x}_t})} \|}^2]}}\nonumber\\
\le& \frac{{4(1 - \beta )(f({{\bar x}_0}) - {f^ * })}}{{\eta d{\sigma _1}T}} + \frac{{\gamma _g^2\sigma _2^6L_{{f_2}}^2}}{{\sigma _1^2}}+ \frac{{4{\kappa _4}{L_{{f_2}}}\eta }}{{n(1 - \beta )d{\sigma _1}}} \nonumber\\
&+\frac{{32{\ell _3}\eta {\kappa _4}}}{{({\rm{1}} - {\rho _1})(1 - \beta )d{\sigma _1}}}\nonumber\\
 =&\frac{{4(1 - \beta )(f({{\bar x}_0}) - {f^*})}}{{\eta d{\sigma _1}T}} + \frac{{\gamma _g^2\sigma _2^6L_{{f_2}}^2}}{{\sigma _1^2}} +\frac{{\eta {d_1}}}{{1 - \beta }}+\frac{\eta {d_2}}{{\gamma _g^2{(1 - \beta) }}},
\end{align}
where
\begin{align}
    {d_1}&={\frac{{16{L_{{f_2}}}d\sigma _2^2L_{{f_1}}^2\sigma _2^2}}{{{\sigma _1}}} + \frac{{128{\ell _3}nd\sigma _2^4L_{{f_1}}^2}}{{({\rm{1}} - {\rho _1}){\sigma _1}}}},\nonumber\\
   {d_2}&= {\frac{{32{L_{{f_2}}}d\sigma _2^2({\gamma _1} + {\vartheta _1})}}{{{\sigma _1}}} + \frac{{256{\ell _3}nd\sigma _2^2({\gamma _1} + {\vartheta _1})}}{{({\rm{1}} - {\rho _1}){\sigma _1}}}}.
\end{align}
By choosing $\eta=(1 - \beta )\sqrt {\frac{n}{T}}$ and running the algorithm for $T>{\tilde \varepsilon _1}$ iterations.
Then, (\ref{888}) can be further simplified as
\begin{align}
\label{section3.2-theorem1-gradient square norm111}
\frac{1}{T}&\sum\limits_{t = 0}^{T - 1} \mathbb{E}{{{ [\| {\nabla f({{\bar x}_t})} \|}^2]}} \nonumber\\
\le& \frac{{4(f({{\bar x}_0}) - {f^*})}}{{\sqrt {nT} d{\sigma _1}}} +  \frac{{\gamma _g^2\sigma _2^6L_{{f_2}}^2}}{{\sigma _1^2}}+ d_1\sqrt {\frac{n}{T}}
+{\frac{d_2}{\gamma _g^2}}\sqrt {\frac{n}{T}}.
\end{align}
This completes the proof. 
\subsection{Proof of Corollary 1}\label{1234567}
When ${\gamma _g} ={T^{ - \frac{1}{8}}}$ and $T > {\tilde \varepsilon _5}$, one can obtain (\ref{section3.2-theorem1-consensus error corollary})--(\ref{section3.2-theorem1-corollary1-gradient square norm}), where 
\begin{align}
{{\tilde c}_1} &={\left(\frac{{640{\phi _1}{{\tilde \varepsilon }_0}}}{\omega }\right )^{\frac{4}{3}}},~~~~~~~
 {\tilde c_2}={\left(\frac{{2880{\phi _1}}}{\omega }\right)^{\frac{4}{3}}},~\nonumber\\
 {\tilde c_3}&={\left(\frac{{54000{\phi _1}{\lambda ^2}}}{{{\omega ^2}{\delta ^2}}}\right)^{\frac{4}{3}}},~~~~
 {\tilde c_4} ={\left(\frac{{48{\phi _1}}}{{{\delta ^3}}}\right)^{\frac{4}{3}}},~\nonumber
\\
 {\tilde c_5}& ={\left(\frac{{2880{\phi _1}{\lambda ^2}}}{{{\delta ^2}\omega }}\right)^{\frac{4}{3}}},~~~~~
 {\tilde c_6}={\left({{128\tilde \varepsilon _0^2{\phi _1}}}\right)^{\frac{4}{3}}},\nonumber\\
 {\tilde c_7}&={\left(\frac{{{5324000}\phi_1}}{{\tilde \varepsilon _8^2{\delta ^6}}}\right)^{\frac{4}{3}}},~~~
  {\tilde c_8}={\left(\frac{{60000{\phi _1}{\lambda ^4}}}{{\omega {\delta ^4}}}\right)^{\frac{4}{3}}},~\nonumber\\
  {\tilde c_9}&={\left(\frac{{40000{\phi _1}}}{{{3\omega ^2}}}\right)^{\frac{4}{3}}},~~~~~
  {\tilde c_{10}}={(8{L_{{f_2}}}{d^2}\sigma _2^2L_{{f_1}}^2)^{\frac{4}{3}}},\nonumber
\end{align}
with $\tilde \varepsilon_0, \phi_1$, $\tilde \varepsilon_8$ and $\tilde \varepsilon_9$ can be found in Appendix~\ref{appendix6.3-Proof of Theorem}.

This completes the proof.
\end{document}